\documentclass{article}

\usepackage{amsfonts,latexsym}

\begin{document}

\newcommand{\h}{{\mathfrak h}}

\newcommand{\D}{\mathcal {D}}

\newcommand{\g}{{\mathfrak g}}

\newcommand{\hgt}{\mbox{ ht }}

\newcommand{\supp}{\mbox{ supp }}

\newcommand{\cst}{{\mathfrak t}}

\newcommand{\mcl}{{\mathfrak l}}

\newcommand{\mfl}{{\mathfrak l}}

\newcommand{\cT}{{\mathcal T}}

\newcommand{\adl}{\mbox{ adj }{\mathfrak l}_i}

\newcommand{\ds}{\displaystyle}

\newcommand{\F}{\frac}

\newcommand{\gzero}{\stackrel{0}{\g}}

\newcommand{\adg}{\mbox{ adj }\gk}

\newcommand{\gk}{\stackrel{k}{\g}}

\newcommand{\Wk}{\stackrel{k}{W}}

\newcommand{\hk}{\stackrel{k}{\h}}

\newcommand{\hkd}{\stackrel{k}{\h^*}}

\newcommand{\Deltak}{\stackrel{k}{\Delta}}

\newcommand{\Pik}{\stackrel{k}{\Pi}}

\newcommand{\alphahat}{\alpha^\vee}
\newcommand{\betahat}{\beta^\vee}

\newcommand{\ahat}{a^\vee}
\newcommand{\ghat}{\g^\vee}
\newcommand{\hhat}{\h^\vee}
\newcommand{\Dhat}{\D^\vee}

\newcommand{\Dk}{\stackrel{k}{\mathcal D}}

\newcommand{\Do}{\stackrel{0}{\mathcal D}}

\newcommand{\ho}{\stackrel{0}{\h}}
\newcommand{\hod}{\stackrel{0}{\h^*}}

\newcommand{\ba}{\begin{array}}
\newcommand{\ea}{\end{array}}

\newcommand{\mbz}{\mathbb Z}
\newcommand{\mbc}{\mathbb C}

 \newtheorem{thm}{Theorem}
\newtheorem{cor}{Corollary}
\newtheorem{lem}{Lemma}
\newtheorem{prop}{Proposition}
\newtheorem{defn}{Definition}
\newtheorem{rem}{Remark}

\input{psfig.sty}

\title{Constructing Graded Lie Algebras}

\author{Meighan I. Dillon\footnote{Partial support was through the Faculty Development Program at Georgia Tech.  Joseph Landsberg's considerable help is especially acknowledged.} }
\date{\today}
\maketitle
\begin{abstract}
The $\mbz$-grading determined by a long simple root of an affine or finite type Lie algebra arises from an adjoint or cominuscule representation of a lower rank semi-simple complex Lie algebra.  Analysis of the relationship between the grading and the representation leads to an extension of Kac's construction of nontwisted affine Lie algebras.\end{abstract}

\begin{flushleft}
\end{flushleft}

\section{Introduction} 
Kac's construction of nontwisted affine Lie algebras produces a $\mbz$-graded infinite dimensional Lie algebra $\g$ from the adjoint representation of a simple complex Lie algebra.  The grading is determined by a so-called special root of $\g$.  We describe two related constructions that produce all $\mbz$-graded finite type and affine Lie algebras--- twisted as well as nontwisted--- where the grading is determined by a long simple root that is not special (an {\em lsn root.}) Each construction starts with a generalized cominuscule representation of a semisimple Lie algebra.  

\begin{defn} \label{com} If $\g=\mcl_1\oplus\cdots\oplus\mcl_t$, where $\mcl_i$ are simple complex Lie algebras, $V=U_1\otimes\cdots\otimes U_t$ is a {\em generalized cominuscule representation} of $\g$ provided $U_i$ is irreducible over $\mcl_i$ with highest weight $n_i\Lambda_i$, $n_i\in{\mathbb Z}_+$, $\Lambda_i$ the fundamental weight associated to a cominuscule simple root of $\mcl_i$. \end{defn}
 
Fix $\gk$, a semi-simple Lie algebra over ${\mathbb C}$.  Let $S$ be a set of positive integers.  For $i\in S$, let  $V_i$ designate an irreducible highest weight representation of $\gk$ with highest weight $\lambda_i$.  The adjoint representation of $\gk$ is $\adg$.
Let $\cT_1$ be the submodule of $\Lambda^2 V_1$ with highest weight(s) of the form $2\lambda_1-\alpha$, $\alpha$ a positive root of $\gk$.  $\cT_1^c$ is the complement of $\cT_1$ in $\Lambda^2 V_1$.  For $i>1$, $\cT_i$ is the submodule of $V_1\otimes V_i$ associated to highest weights greater than or equal to those of the form $\lambda_1+\lambda_i-\alpha$, where $\alpha$ is any  positive root of $\gk$.   $\cT_i^c$ is the complement of $\cT_i$ in $V_1\otimes V_i$.

\begin{thm}\label{main1}{{\bf The Affine Algorithm}}
Suppose $V=V_1$ is a generalized cominuscule representation of $\gk$. \begin{enumerate}  Set $j=1$.
\item If $\cT_j^c$ is irreducible, let $V_{j+1}=\cT_j^c$.  Increase $j$ by $1$ and repeat this step.
\item If $\cT_j^c=\adg\oplus\;{\mathbb C}$, let $V_{j+1}=\adg$.  Increase $j$ by $1$ and go to step (4).
\item Otherwise $V$ is inadmissible.

\item For $i>j$, let $V_i=\adg$ if $i\equiv 0$ mod $j$.  Otherwise take $\ell\in\{1,\ldots,j-1\}$ so that $\ell\equiv i$ mod $j$ and let $V_i\cong V_\ell$.  For $-i<0$, take $V_{-i}=V_i^*$.   \end{enumerate}

The algorithm is effective, terminating with $j\leq 6$. 

If $V$ is admissible, let $\ds\g_0=\gk\oplus\;{\mathbb C}\oplus {\mathbb C}$ and, for $i\neq 0$, let $\g_i=V_i$ .  $\g=\bigoplus_{i\in{\mathbb Z}}\g_i$ is a ${\mathbb Z}$-graded affine algebra, the grading determined by an lsn root of $\g$.  Moreover, any lsn-graded affine Lie algebra can be constructed this way.   
\end{thm}

The next result generalizes the Minuscule Algorithm as detailed in \cite{lan}.  

\begin{thm}\label{main2}{{\bf The Finite Algorithm}} 
Suppose $V=V_1$ is a generalized cominuscule representation of $\gk$.  
\begin{enumerate}  Set $j=1$.
\item If $\cT_j^c=\{0\}$, go to step (4).
\item If $\cT_j^c$ is irreducible, take $V_{j+1}=\cT_j^c$.  Increase $j$ by $1$ and go back to step (1).  
\item Otherwise $V$ is inadmissible.
\item For $i\in\{-1,\ldots,-j\}$, take $V_{-i}=V^*_i$. For $i>j$ and $i<-j$, take $V_i=\{0\}$. \end{enumerate}

The algorithm is effective, terminating with $j\leq 6$.

If $V$ is admissible, let $\ds\g_0=\gk\oplus\;{\mathbb C}$ and, for $i\neq 0$, let $\g_i=V_i$.  
$\g=\bigoplus_{i\in{\mathbb Z}}\g_i$ is a ${\mathbb Z}$-graded finite type Lie algebra, the grading determined by an lsn root of $\g$.  Moreover, any lsn-graded finite type Lie algebra can be constructed this way.   
\end{thm}

\begin{rem} Both algorithms relate to work of B. Kostant.  Let $\g_0$ be a reductive Lie algebra with module ${\mathfrak p}$.  Theorem~1.50 in \cite{kos} gives necessary and sufficient conditions for determining whether a Lie bracket can be imposed on $\g=\g_0+{\mathfrak p}$, consonant with the action of $\g_0$ on $\mathfrak p$.  Our work starts with
semi-simple $\gk$ and a representation $V$ and constructs a minimal $\mbz$-graded representation space $\mathfrak p=\oplus\g_i$ so that $\g_{-i}=\g_i^*$, $\g_1=V$, $\g_0$ contains $\gk$, and $\g_0+{\mathfrak p}$ is an affine or finite type Lie algebra with an lsn-grading.
\end{rem}

The rest of the paper is dedicated to background analysis and the proofs of the theorems.

\section{Terminology and Notation}
A Lie algebra of {\em finite type} is a simple, finite dimensional Lie algebra over ${\mathbb C}$. A {\em semi-simple Lie algebra} is always complex semi-simple.   All Lie algebras here are affine or semi-simple.  Given a Lie algebra $\g$, we fix a Cartan subalgebra $\h$. $\Delta$ is the set of roots of $\g$, $\Delta_+$ the positive roots, $\Pi$ the set of simple roots, etc.  The Dynkin diagram for $\g$ is $\D$.  We number the simple roots of a finite type Lie algebra as in \cite{bou} and \cite{hum}.  Otherwise, our notation and conventions typically follow \cite{kac}.  Unless indicated to the contrary, $\alpha_i$ is a simple root.  The  extra root on an extended (or affine) Dynkin diagram (cf. \cite{bou}) is always $\alpha_0$.  The {\em labels} of $\D$ are as in the Tables Aff in Chapter~4 of \cite{kac}.  It is often convenient to identify elements of $\Pi$ with the nodes of $\D$. 

The ${\mathbb Z}$-span of the roots of a Lie algebra form a lattice designated $Q$. ${\mathbb Z}_+$-linear combinations of simple roots comprise $Q_+$.  If $\beta$ is in $Q_+$, write it as a ${\mathbb Z}_+$-linear combination of simple roots. The {\em support} of $\beta$ is then the set of simple roots with positive coefficients.  The {\em height} of $\beta$ is the sum of those coefficients.  When applied to roots, the words {\em highest} and {\em lowest} refer to height. Let $a_i$ be the label of $\D$ associated to $\alpha_i$. In the finite case, the highest root of $\g$ is $\theta=\sum_{i=1}^\ell a_i\alpha_i$; in the affine case, the minimal positive imaginary root is $\delta=\sum_{i=0}^\ell a_i\alpha_i$. 

A simple root of $\g$ is {\em special} if it is conjugate to $\alpha_0$ under a diagram automorphism. This makes sense for a finite type $\g$ if we consider the extended Dynkin diagram for $\g$.  A {\em cominuscule root} is a special root for a finite type Lie algebra (cf. \cite{mag}.)  The cominuscule roots are:  all $\alpha_i$ in $A_n$; $\alpha_1$ in $B_n$; $\alpha_n$ in $C_n$; $\alpha_1$, $\alpha_{n-1}$ and $\alpha_n$ in $D_n$; $\alpha_1$ and $\alpha_6$ in $E_6$; $\alpha_7$ in $E_7$. 

Let $V$ be a representation of $\g$.  A weight vector $v^+$ is a {\em highest weight vector} of $V$ provided $e_\alpha.v^+=0$ for all positive root vectors $e_\alpha$ in $\g$.  The weight associated to $v^+$ is then a {\em highest weight}.   $f_\alpha$ designates a root vector in $\g$ associated to $-\alpha$, where $\alpha$ is a positive root.  A weight vector $v^-$ is a {\em lowest weight vector} provided $f_\alpha.v^-=0$ for all positive roots $\alpha$.  To indicate that $V$ is an irreducible $\g$-module with highest weight $\lambda$, we sometimes write $V(\lambda)$ instead of $V$. 

Suppose $\g$ is finite type or affine.  Distinguishing a long simple root of $\g$, $\alpha_k$, we get a ${\mathbb Z}$-grading on $\g$ by $\deg e_j=-\deg f_j=\delta_{jk}$ (Kronecker delta).  (See \cite{kac}, \S 1.5.) Here $e_i$, $f_i$ are root vectors of $\g$ associated to simple $\alpha_i$ and $-\alpha_i$ respectively.  The $\mbz$-grading induced by $\alpha_k$ is the {\em $\alpha_k$-grading of $\g$}.  If $\alpha_k$ is an lsn root, we say the $\alpha_k$-grading  of $\g$ an {\em lsn-grading.}  

Mark the nodes of $\D$ adjacent to $\alpha_k$. If a marked  node represents a shorter root, label it with the number of edges between it and $\alpha_k$.  Excising $\alpha_k$ and adjacent edges from $\D$ we get $\Dk$, a marked and labeled Dynkin diagram associated to a semi-simple Lie algebra $\gk$ (cf. \cite{oni}, Chapter 3, \S3.5.) 
When $\g$ is affine, $\gzero$ is the {\em underlying algebra of finite type.} 
$\Dk$ represents an irreducible module over $\gk$.  Its highest weight is the sum of the fundamental weights associated to the marked nodes, with multiplicity according to labels. Note that in the nontwisted affine algebras, $\Do$ is the diagram for the adjoint representation of $\gzero$.  In general the module associated to $\Dk$ corresponds to $\g_{-1}$, the sum of the root spaces of $\g$ associated to roots with coefficient $-1$ on $\alpha_k$. The connection between the $\alpha_k$-grading and the representation it determines is discussed thoroughly in \S4 below.

Fix $\hk\subset\h$, a Cartan subalgebra for $\gk$. The set of roots of $\gk$ is $\Deltak$, the set of simple roots is $\Pik$, etc. Fix $\Pik$ so there is a one-to-one correspondence between $\Pik$ and $\Pi-\{\alpha_k\}$ and use the same symbols to designate $\alpha_j$ in $\Pik$ and in $\Pi$, for $j\neq k$.  
 Let $(.,.)$ designate a fixed standard invariant form on $\g$.  Sometimes we call it the Killing form. Normalize the form so that $\|\alpha\|^2=2$ for any long root $\alpha$ of $\g$.  We use the same notation for the form as it restricts to $\h$, also to $\gk$ and $\hk$. The canonical isomorphism determined by the form is $\nu:\h\rightarrow\h^*$ and $(.,.)$ designates the induced form on $\h^*$ as well.    Let $<.,.>$ denote the pairing of $\g$ and $\g^*$, its dual.  If $\alpha\in\Delta$, $\alphahat\in\h$ designates its coroot.  If $\alpha_i$ is long, $-<\alpha_i,\alphahat_j>$ is the number of edges between $\alpha_i$ and $\alpha_j$ in $\D$.  For simple roots $\alpha_i$, $\alpha_j$, $\ds <\alpha_i,\alphahat_j>=\F{2(\alpha_i,\alpha_j)}{\|\alpha_j\|^2}.$  Also $\ds\nu(\alphahat_i)=\F{2\alpha_i}{\|\alpha_i\|^2}.$ 

\section{Brackets}
Here we consider brackets on the $\alpha_k$-graded pieces of $\g=\oplus_i\g_i$, looking at them in terms of the action of $\g_0$ on $\g_i$ by $x.g=[x\;g]$.  It is enough to consider $[\g_1\;\g_j]$ and $[\g_{-1}\;\g_j]$ as other brackets are defined iteratively in terms of these. 

If $\g$ is affine type, $\g_0=\gk\oplus\;{\mathbb C}d_k + {\mathbb C}K$, where ${\mathbb C}K$ is the center of $\g$ and $d_k=\nu^{-1}(\Lambda_k)$, $\Lambda_k$ the fundamental weight associated to $\alpha_k$.  $\Lambda_k$ and $d_k$ are isotropic and $(K,d_k)=a_k$, $a_k$ the label on $\D$ associated to $\alpha_k$.  Note that $K$ belongs to the derived algebra of $\g$ but $d_k$ does not. (See \cite{kac}.)
If $\g$ is finite type, $\g_0=\gk\oplus\;{\mathbb C}d_k$, where $d_k\in\h$ is given by $\nu(d_k)=\Lambda_k=\sum_it_i\alpha_i$, for scalars $t_i\in\;{\mathbb C}$.  If $x\in\g_t$, $[d_k\;x]=tx$ and $(\Lambda_k,\Lambda_k)=(d_k,d_k)=t_k$.  

\subsection{{\mathversion{bold}
$[\cdot\;,\;\cdot]:\g_{-1}\otimes\g_1\rightarrow\g_0$}}

  Let $\{X_i\}$ and $\{Y_i\}$ be Killing dual bases of $\gk$. If $\g$ is affine, $\{X_i\}\cup\{\F{1}{a_k}d_k, K\}$ and $\{Y_i\}\cup\{K,\F{1}{a_k}d_k\}$ are dual bases of $\g_0$. For $u_{-1}\in\g_{-1}$, $u_1\in\g_1$, we have
$$[u_{-1}\;u_1]=\sum_i (X_i,[u_{-1}\;u_1])Y_i+\F{1}{a_k}(d_k,[u_{-1}\;u_1])K+\F{1}{a_k}(K,[u_{-1}\;u_1])d_k=$$$$-\sum_i(X_i,[u_1\;u_{-1}])Y_i+\F{1}{a_k}([d_k\;u_{-1}],u_1)K+\F{1}{a_k}([K\;u_{-1}],u_1)d_k=$$$$-\sum_i([X_i\;u_1],u_{-1})Y_i-\F{1}{a_k}(u_{-1},u_1)K.$$ Then in the affine case we have, $[u_{-1}\;u]=-\sum_i(u_{-1},X_i.u_1)Y_i-\F{1}{a_k}(u_{-1},u_1)K.$

If $\g$ is finite type, $\{X_i\}\cup\{\F{1}{t_k}d_k\}$ and $\{Y_i\}\cup\{d_k\}$ are dual bases of $\g_0$.  For $u_{-1}\in\g_{-1}$, $u_1\;\in\;\g_1$, a calculation similar to the one above yields $[u_{-1}\;u_1]=-\sum_i(u_{-1},X_i.u_1)Y_i-\F{1}{t_k}(u_{-1},u_1)d_k.$

\subsection{{\mathversion{bold}
$[\cdot\;,\;\cdot]:\Lambda^2\g_{1}\rightarrow\g_2$}}

Let $u$, $v\;\in\;\g_1$, $u_{-1}\;\in\;\g_{-1}$.  In the affine case, we have
$$[u_{-1}\;[u\;v]\;]=[\;[u_{-1}\;u]\;v]+[u\;[u_{-1}\;v]\;]=[\;[u_{-1}\;u]\;v]-[\;[u_{-1}\;v]\;u]=$$$$
-\sum_i(u_{-1},X_i.u)Y_i.v-\F{1}{a_k}(u_{-1},u)K.v+\sum_i(u_{-1},X_i.v)Y_i.u$$$$+\F{1}{a_k}(u_{-1},v)K.u=\sum_i(u_{-1},X_i.v)Y_i.u -\sum_i(u_{-1},X_i.u)Y_i.v.$$  Viewing $[u\;v]$ as an element in $\Lambda^2\g_1$, we then have \begin{equation}\label{3}[u\;v]=\sum_i X_i.v\wedge Y_i.u.\end{equation}
A similar calculation in the finite case gives us 
$[u_{-1}\;[u\;v]\;]=$
$$-\sum_i(u_{-1},X_i.u)Y_i.v-\F{1}{t_k}(u_{-1},u)d_k.v+\sum_i(u_{-1},X_i.v)Y_i.u+\F{1}{t_k}(u_{-1},v)d_k.u$$$$=\sum_i(u_{-1},X_i.v)Y_i.u -\sum_i(u_{-1},X_i.u)Y_i.v+\F{1}{t_k}((u_{-1},v)u-(u_{-1},u)v).$$
\begin{equation}\label{4}\mbox{Then in the finite case   }[u\;v]=\sum_i X_i.v\wedge Y_i.u+\F{1}{t_k}v\wedge u.\end{equation}

\subsection{{\mathversion{bold}
$[\cdot\;,\;\cdot]:\g_{-1}\otimes\g_t\rightarrow\g_{t-1}$} and {\mathversion{bold}
$[\cdot\;,\;\cdot]:\g_{1}\otimes\g_t\rightarrow\g_{t+1}$}}

Let $\g$ be affine type and suppose $t>1$.  Assume $[\g_{-1}\;\g_{t}]$ is given by 
$$[u_{-1}\;u_1\;\ldots\;u_t]=-\sum_i(u_{-1},X_i.u_1)Y_i.[u_2\;\ldots\;u_t]-(u_{-1},X_i.u_2)[u_1\;Y_i.[u_3\;\ldots\;u_t]\;]-$$$$\ldots-(u_{-1},X_i.u_{t-1})[u_1\;\ldots\;u_{t-2}\;Y_i.u_t]+(u_{-1},X_i.u_t)[u_1\;\ldots\;Y_i.u_{t-1}].$$  At the same time, assume that $[\g_1\;\g_{t}]$ is given by 
$$[u_1\;\ldots\;u_t]=-\sum_iX_i.u_1\otimes Y_i.[u_2\;\ldots\;u_t]- X_i.u_2\otimes [u_1\;Y_i.[u_3\;\ldots \;u_t]\;]-$$$$\ldots-X_i.u_{t-1}\otimes [u_1\;\ldots\;u_{t-2}\;Y_i.u_t]+X_i.u_t\otimes [u_1\;\ldots \;u_{t-2}\; Y_i.u_{t-1}].$$
By induction and the Jacobi identity we have $$[u_{-1}\;u\;u_1\ldots\;u_t]=[\;[u_{-1}\;u]u_1\;\ldots\;u_t]+[u\;u_{-1}\;u_1\;\ldots\;u_t]=$$$$-\sum_i(u_{-1},X_i.u)\otimes Y_i.[u_1\;\ldots\;u_t]-\sum_i(u_{-1},X_i.u_1)[u\;Y_i.[u_2\;\ldots\;u_t]\;]-$$$$\ldots-(u_{-1},X_i.u_{t-1})[u\;u_1\;\ldots\;u_{t-2}\;Y_i.u_t]+(u_{-1},X_i.u_t)[u\;u_1\;\ldots\;Y_i.u_{t-1}].$$ From there it follows that 
$$[u\;u_1\;\ldots\;u_t]=-\sum_iX_i.u\otimes Y_i.[u_1\;u_2\;\ldots\;u_t]- X_i.u_1\otimes [u\;Y_i.[u_2\;\ldots \;u_t]\;]-$$\begin{equation}\label{5}\ldots-X_i.u_{t-1}\otimes [u\;u_1\;\ldots\;u_{t-2}\;Y_i.u_t]+X_i.u_t\otimes [u\;u_1\;\ldots \;u_{t-2}\; Y_i.u_{t-1}].\end{equation}

Let $\g$ be finite type and suppose $t>1$.  Applying induction and the Jacobi identity as above to both $[u_{-1}\;u_1\;\ldots\;u_t]$ and $[u_1\;\ldots\;u_t]$ we get
$$[u\;u_1\;\ldots\;u_t]=-\sum_iX_i.u\otimes Y_i.[u_1\;u_2\;\ldots\;u_t]-X_i.u_1\otimes [u\;Y_i.[u_2\;\ldots \;u_t]\;]-$$$$\ldots-X_i.u_{t-1}\otimes [u\;u_1\;\ldots\;u_{t-2}\;Y_i.u_t]+X_i.u_t\otimes [u\;u_1\;\ldots \;u_{t-2}\; Y_i.u_{t-1}]$$$$-\F{t}{t_k}u\otimes[u_1\;\ldots\;u_t]-\F{t-1}{t_k}u_1\otimes [u\;u_2\;\ldots\;u_t]-$$\begin{equation}\label{6}\ldots-\F{1}{t_k}u_{t-1}\otimes [u\;u_1\;\ldots\;u_t]+\F{1}{t_k}u_t\otimes[u\;u_1\;\ldots\;u_{t-1}].\end{equation}

\section{The $\alpha_k$-Grading}
This section is an elaboration on ideas sketched in \cite{oni}, Chapter~3, \S3.5.

The $i$th $\alpha_k$-graded piece of $\g$, $\g_i$, is a sum of root spaces of $\g$ associated to roots $\alpha$ with the form 
$\ds\alpha=i\alpha_k +\sum_{j\neq k} c_j \alpha_j.$  
Note that $\g_{\alpha}\subset \g_i$ if and only if $\g_{-\alpha}\subset\g_{-i}$. 
If $\g$ is affine, suppose it is of type $X_N^{(r)}$.  Recall that $\alpha$ is a real root of $X_N^{(r)}$ if and only if $\alpha+nr\delta$ is a real root for all $n$ in ${\mathbb Z}$.
In the finite case, $a_k$ is maximal so that $\g_{a_k}$ is nonzero.   In the affine case, $\g_{ra_k}$ contains $\g_{r\delta}$, along with some positive real root spaces.  The lowest root associated to $\g_{ra_k}$ has the form $r\delta-\theta$, where $\theta$ is a highest root of $\gk$. ($\gk$ is semi-simple in general so there may be more than one highest root.) This is an immediate consequence of the definition of the $\alpha_k$-grading along with Proposition~6.3~(d) in \cite{kac}. 
If $\g$ is finite type, let $ra_k=a_k$.

Let $\sum_{i=0}^\ell c_i\alpha_i=\alpha\in\Delta_+$.  The {\em root diagram} ${\D}(\alpha)$ is the labeled subdiagram of $\D$ comprised of the nodes and connecting edges associated to simple roots in the support of $\alpha$. If $c_i>1$, the node associated to $\alpha_i$ is labeled $c_i$. We can extend the notion of a root diagram to apply to any element of $Q_+$. We distinguish $\alpha_k$ in ${\D}(\alpha)$ by coloring the associated node. 
Say $\beta\;\in\;\Delta_+$ is a {\em subroot} of $\alpha\;\in\;\Delta_+$ provided $\alpha-\beta$ is in $Q_+$. 

\begin{lem} \label{subrootlemma} Let $\gamma_i$ be a lowest root associated to $\g_{i}$ so that $\gamma_1=\alpha_k$.  

$$\gamma_2-(\gamma_1+\alpha_k)=\gamma_2-2\alpha_k=\sum_{j\neq k}c_j\alpha_j=\beta_1+\beta_2+\beta_3$$ where $\beta_i\in\;\Deltak_+$ have the following properties.

\begin{enumerate}
\item $\beta_1=\sum_{j\neq k}c'_j\alpha_j$ is a maximal subroot of $\sum_{j\neq k}c_j\alpha_j$, in the sense that $\beta_1$ is a root, $0\leq c'_j\leq c_j$ for all $j$, and if $c'_j$ were replaced with $c'_j+1$, for any $j$, then $\beta_1$ would no longer be a root, or $c'_j+1$ would exceed $c_j$.  $\beta_2$ is a maximal subroot of $\sum_{j\neq k}(c_j-c'_j)\alpha_j\in Q_+$ and $\beta_3=\sum_{j\neq k}c_j\alpha_j-\beta_1-\beta_2$.

\item If $i<j$ and $\beta_i$ and $\beta_j$ come from the same simple component of $\gk$, then ${\D}(\beta_j)\subseteq {\D} (\beta_i)$.

\item  $S=\{\beta_1,\beta_2,\beta_3\}$ is an {\em inert set of roots}, that is, neither the sum nor the difference of elements in $S$ is a root.  In particular, $(\beta_i,\beta_j)=0$ for $i\neq j$.

\item  For any arrangement of indices, $\gamma_1+\beta_i$, $\gamma_1+\beta_i+\beta_j$  and $\gamma_1+\beta_i+\beta_j+\beta_k$ are all roots of $\gk$. 

\item $\beta_3$ is simple.
 \end{enumerate}

If $2<i\leq ra_k$, then $\gamma_i-(\gamma_{i-1}+\alpha_k)=\beta_1+\beta_2$ where $\beta_1, \beta_2\in\Deltak_+$.  {\em Mutatis mutandis}, $\gamma_i$, $\beta_1$, and $\beta_2$ satisfy the first four properties above.
\end{lem}

{\bf Proof}  Any element of $Q_+$ is a sum of successively maximal subroots as per the first statement of the lemma.  We have to show that there are  three such subroots for$\gamma_2$ and, when $2<i\leq ra_k$, two for $\gamma_i$.  We do so by demonstrating how to construct the $\beta_i$s in each case.   It is a routine matter to check that the other items in the lemma then follow.  Constructing the $\beta_i$s is straightforward once we identify the possible root diagrams for $\gamma_i$.

The tables of roots in \cite{bou} give us $\gamma_i$ in the finite cases.  In the affine cases, we appeal to Proposition~6.3 in \cite{kac}. If $i<a_k$, $\gamma_i$ is a finite type root, that is, its support is contained in $\Deltak$ and it appears on a list in \cite{bou}.  If $i=ja_k$, $\gamma_i=j\delta-\theta$, where $\theta$ is a highest long or short root of $\gk$. (The root is short if $j<r$.) The one case left to consider is when $i>a_k$ is not a multiple of $a_k$.  That happens only if $\g=E_6^{(2)}$, $a_k=2$, and $i=3$.  Then $\gamma_3=j\delta-\xi$, where $\xi$ is the highest short root of $\g$ with coefficient $1$ on $\alpha_k$. It turns out this is also a root of finite type.  It appears as Case~5 in Figure~2. 

We list all possible root diagrams for $\gamma_2$ in Figure~1. Along with each root diagram ${\D}(\gamma_2)$ is the decomposition of $\gamma_2-2\alpha_k$ into a sum of three maximal subroots.  These disconnected root diagrams are ${\D}(\beta_1)$, ${\D}(\beta_2)$, ${\D}(\beta_3)$ in each case. 

Figure~2 shows the possibilities for ${\D}(\gamma_3)$ when $ra_k\geq 3$. Note that at this stage, all roots are of finite type. 

\psfig{figure=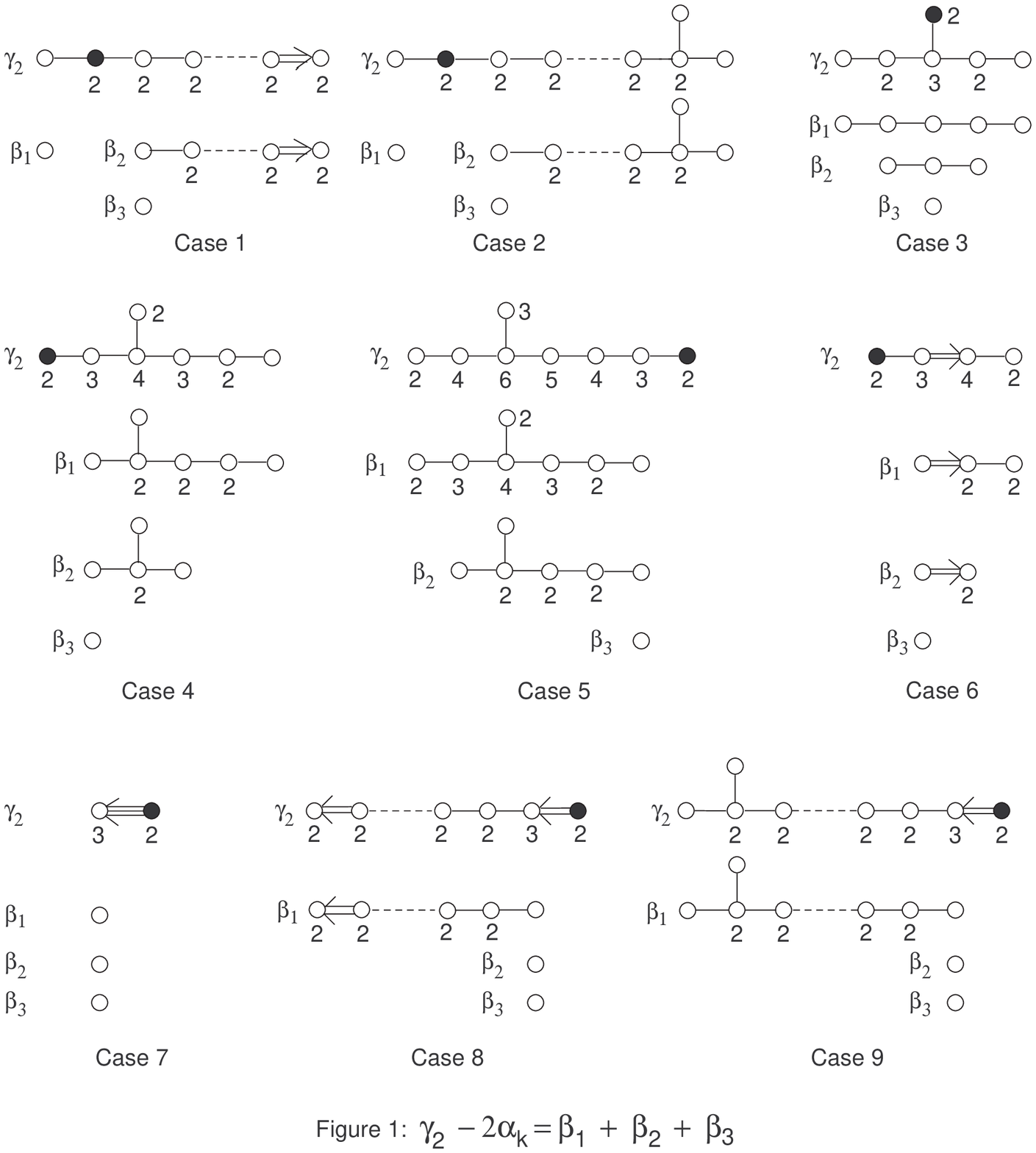,height=4in}

\psfig{figure=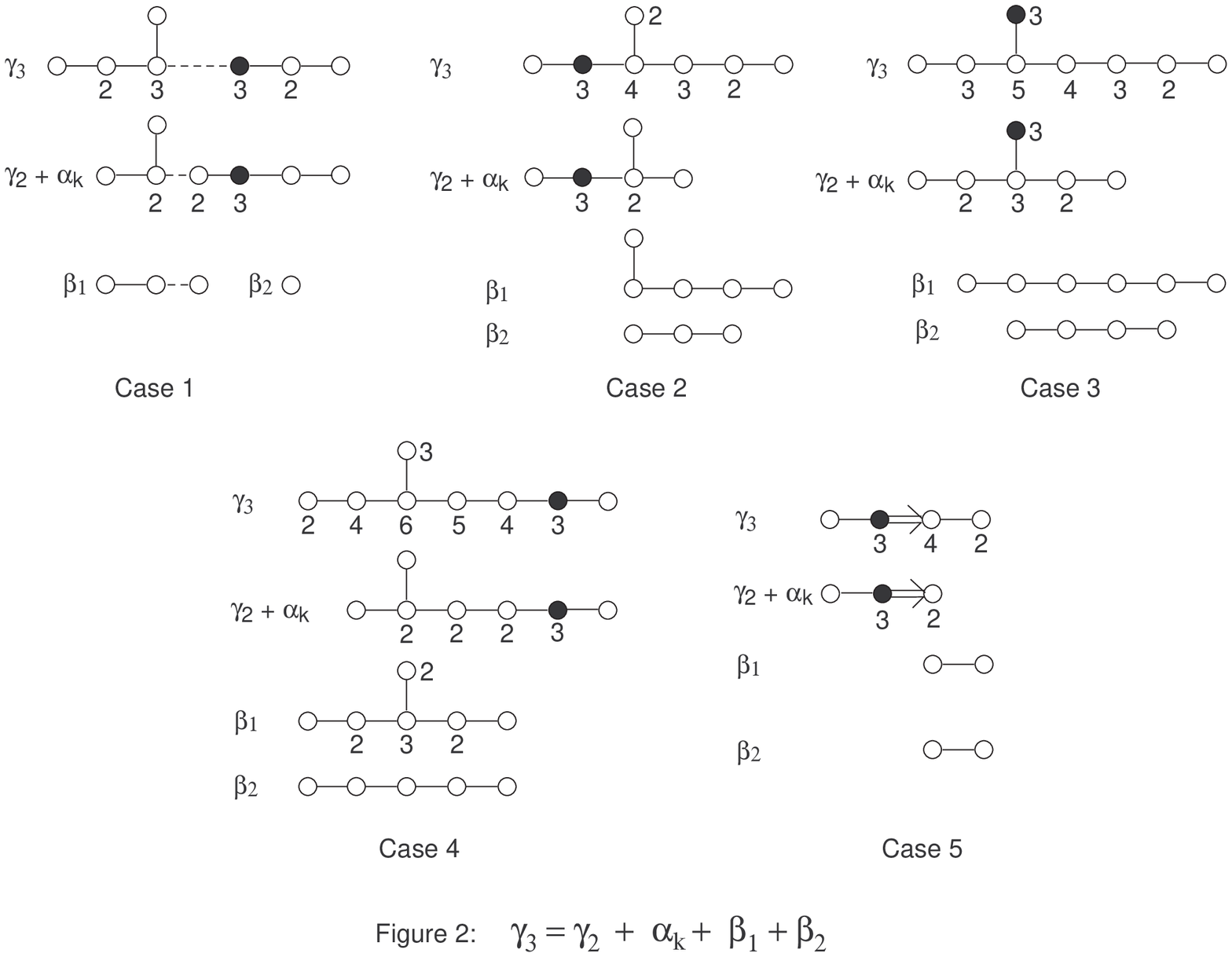,height=3.5in}

 Along with each root diagram is the decomposition of $\gamma_3$ into $\gamma_2+\alpha_k$ and two maximal subroots.  As in Figure~1, these disconnected root diagrams are ${\D}(\beta_1)$ and ${\D}(\beta_2)$ in each case.  In these and higher order cases, we have to match $\gamma_2$ to $\gamma_3$.  With the catalog given in Figure~1, this is not difficult.  

A complete catalog of root diagrams goes through $a_k=6$.  We leave the remaining three figures to the reader. 
$\Box$

\subsection{The Affine Case} When $\g$ is affine, $\nu(K)=\delta$.  Designating by $\ahat_i$ the labels on the diagram for the Lie algebra dual to $\g$, we have $K=\sum_{i=0}^\ell \ahat_i\alphahat_i$:  in general $\nu(\alphahat_i)=\F{a_i}{a^\vee_i}\alpha_i$ (cf. \cite{kac}, \S6.2.)  Note finally that $\h^*=\hkd\oplus\; {\mathbb C}\delta+{\mathbb C}\Lambda_k.$ 

Next is a record of some facts we use throughout our discussion.  All are either stated explicitly in \cite{kac} or easily deduced.

\begin{lem}\label{affinefactslem}
$$\begin{array}{lllr}
\nu(\alphahat_k)=\alpha_k&\nu(K)=\delta&\nu(d_k)=\Lambda_k&(d_k, d_k)=0\\ \\

 (\Lambda_k,\Lambda_k)=0& (K,K)=0& (\delta,\delta)=0&(\alphahat_i,d_k)=\delta_{ik}\\\\
 (\Lambda_k,\alpha_k)=1&(K,d_k)=a_k& (\Lambda_k,\delta)=\ahat_k&\Lambda_k(K)=\ahat_k\\ \\
<\delta,\alphahat_i>=0&<\delta , d_k>=a_k& <\alpha_i,d_k>=\delta_{ik}& \Box
\end{array}$$ 
\end{lem} 

Let $\pi:\h^*\rightarrow \hkd$ be given by
$\pi(\alpha_i)=\alpha_i$, for $i\neq k$, $\pi(\delta)=0$, and $\pi(\Lambda_k)=0$.  Since $\delta=\sum_i a_i\alpha_i$, this gives us
$0=\pi(a_k\alpha_k+\sum_{i\neq k}a_i\alpha_i)=a_k\pi(\alpha_k)+\sum_{i\neq k}a_i\alpha_i$
so that 
$\pi(\alpha_k)=-\sum_{i\neq k}\F{a_i}{a_k}\alpha_i=-\F{1}{a_k}\delta+\alpha_k.$  

Note that $\alpha_k$ is special if and only if $ra_k=1$.  

\begin{lem} \label{afg1ir} $\g_{1}$ is an irreducible $\gk$ module. When $\alpha_k$ is  special, $\g_{1}$ is isomorphic to the adjoint representation of $\gk$.  
\end{lem}
{\bf Proof} $\gk$ acts on $\g_{1}$ via the bracket.  That this defines $\g_{1}$ as a $\gk$ module follows from the Jacobi identity on $\g$ and the observation that $[\g_0\;\g_{1}]\subset \g_{1}$ while $\gk\subset\g_0$. As a $\gk$ module, $\g_1$ has lowest weight $\pi(\alpha_k)=-\F{1}{a_k}\delta+\alpha_k$.  
 We show that up to scalar multiples, $e_{\alpha_k}=e_k$ in $\g_{\alpha_k}$ is the unique lowest weight vector for $\g_{1}$.  From there, it follows that $\g_{1}$ is irreducible.
Assume first that $ra_k>1$. 

Suppose $\pi(\beta)$ is another lowest weight associated to $\g_1$ so that for $\beta\in\Delta_+$, $e_\beta\in\g_\beta$, $f_j.e_\beta=0$ for all $\alpha_j\in\Pik$. We claim there must be $e_i\in\gk$ with $e_i.e_\beta\neq 0$.  

By Lemma~2.1 in \cite{dil}, we can write $\beta=\alpha_k+\sum_{j=1}^m\alpha_{i_j}$ where $\alpha_{i_j}\in\;\Pik$ and for all $t\in\{1,\ldots,m\}$, $\alpha_k+\sum_{j=1}^t\alpha_{i_j}\;\in\;\Delta_+.$ We also know that if $\beta$ and $\beta-\alpha_j$ are roots, then $f_j.\g_{\beta}\neq 0$ (cf. \cite{kac}, Prop.~3.6.)  If $f_j.e_\beta=0$, then $\dim \g_\beta>1$, so $\beta$ is an imaginary root and $a_k=1$. Since $ra_k>1$, it follows that $\g$ is twisted affine. Now $\alpha_k$ and $-\alpha_k$ are long roots, and in the twisted affine algebras, $\delta+\alpha$ is not a root when $\alpha$ is long (cf. \cite{kac}, Prop.~6.3.) So if $f_i.e_\beta=0$ for all $\alpha_i\in\Pik$, then $f_k.e_\beta$ is also zero.  By the same argument, if we assume $e_i.e_\beta=0$ for all $\alpha_i\in\Pik$, then $e_k.e_\beta=0$ since $\beta+\alpha_k$ is not a root.  This gives us a root vector, $e_\beta$, that commutes with all root vectors and everything in $\h$ except $d_k$.  This means $e_\beta$ is a multiple of $K$, which is absurd.  Conclusion:  if we insist that $f_i.e_\beta=0$ for all $i\in\Pik$, then (1) $\beta$ is imaginary, and (2) there must be $e_i\in \gk$ with $e_i.e_\beta\neq 0$.

Since $\beta+\alpha_i$ is real, the argument above gives us $f_j$ with $f_j.e_i.e_\beta\neq 0$. If $\alpha_j=\alpha_i$ we have $f_j.e_i.e_\beta=-\alphahat_i.e_\beta+e_i.f_j.e_\beta= e_i.f_j.e_\beta\neq 0$ as $<\beta,\alphahat_i>=0$.  If $\alpha_j\neq\alpha_i$, $f_j$ and $e_i$ commute and we again get $e_i.f_j.e_\beta\neq 0$, contradicting $f_j.e_\beta=0$.  Final conclusion in case $ra_k>1$:  up to scalar multiples, $e_k$ is a unique lowest weight/root vector associated to $\g_1$ so $\g_1$ is irreducible.

Next suppose $ra_k=1$ so that $\alpha_k$ is special and $\g$ is nontwisted affine.  $\gk$ is then simple and $\g_\delta\subset\g_1$.  The real roots of $\g_1$ are precisely those of the form $\alpha+\delta$, where $\alpha$ is any root of $\gk$:  this follows from \cite{kac}, Prop. 6.3.  Thus there is a bijective correspondence between the roots of $\gk$ and the real roots/weights of $\g_1$.  In particular, the highest weight of $\g_1$ is $\pi(\theta+\delta)=\theta$, the highest root of $\gk$.  It follows that there is a copy of $\gk$ inside $\g_1$.  Note that all real root/weight spaces of $\g_1$ are one dimensional and that the dimension of $\g_\delta$ is the rank of $\gk$.  By dimension, $\g_1$ must be isomorphic to the adjoint representation of $\gk$.  (Note that  $\g_\delta$ is a Cartan subalgebra.) Since $\gk$ is simple, it follows that $\g_1$ is irreducible in this case as well.    $\Box$ 

\begin{lem}\label{afdual}
$\g_i$ and $\g_{-i}$ are dual representations of $\gk$.
\end{lem}

{\bf Proof} The action of $\gk$ on $\g_i$ is via the bracket and it follows that $\g_i$ is a $\gk$ module.  As $\g_0$ is the direct sum of $\gk$ and two copies of the trivial representation of $\gk$, it is self-dual. 
Assume $i>0$.  A lowest root associated to $\g_{-i}$ has the form
$-\beta=-i\alpha_k-\sum_{j\neq k} c_j\alpha_j$
where $\sum_{j\neq k} c_j$ is maximal so that $\beta$ is a positive root of $\g$. This is precisely the criterion that determines that $\beta$ is a highest root associated to $\g_i$.  
Thus as $\gk$ modules, $\g_{i}$ has a highest weight $\pi(\beta)$ if and only if $\g_{-i}$ has lowest weight $\pi(-\beta)=-\pi(\beta)$.  This is to say that for $i\neq 0$, $\g_i$ and $\g_{-i}$ are dual $\gk$ modules.   $\Box$ 

\begin{lem} \label{grdel} As a $\gk$ module, $\g_{ra_k}$ is isomorphic to the adjoint representation of $\gk$. 
\end{lem}

{\bf Proof}  Since $\delta=\sum_{i=0}^\ell a_i \alpha_i$, $\g_{r\delta}\subset \g_{ra_k}$.  All other roots $\beta=a_k\alpha_k+\sum_{j\neq k}c_j\alpha_j$ associated to $\g_{ra_k}$ are real in which case $\beta-r\delta$ is a real root of $\g$ and the associated root space is one dimensional.  Since $\alpha_k$ is not in the support of $\beta-r\delta$, $\g_{\beta-r\delta}\subset\gk$. Conversely, if $\g_\beta$ is a root space in $\gk$, then $r\delta-\beta$ is a real root of $\g_{ra_k}$.  Thus root spaces in $\gk$ are in one-to-one correspondence with real root spaces of $\g_{ra_k}$.  In particular, highest weights of $\g_{ra_k}$  correspond precisely to the highest roots of $\gk$. (In general, $\gk$ is semi-simple so there may be more than one highest root.) Since $\delta|_{\hk}=0$ and $\dim \g_{r\delta}=\mbox{ rank }\gk$ (cf. \cite{kac} Corollary~8.3), $\g_{r\delta}\subset\g_{ra_k}$ corresponds to a Cartan subalgebra of $\gk$.  It follows that $\g_{ra_k}$ is isomorphic to the adjoint representation of $\gk$.  $\Box$  

\begin{lem}\label{gigjcon}
For $i,j\neq 0$,  $i\equiv j\mbox{ mod }ra_k$ if and only if, as $\gk$ modules, $\g_i\cong\g_j$.
\end{lem}

{\bf Proof} Assume $i$ and $j$ are nonzero and that $i\equiv j\mbox{ mod }ra_k$.

Lemma~\ref{grdel} and its proof remain valid if we replace $r\delta$ with $nr\delta$ and $ra_k$ with $nra_k$ for $n\in{\mathbb Z}$.  In particular, $\beta\in\h^*$ determines a root space in $\g_{nra_k+i}$ if and only if $\beta-nr\delta$ determines a root space in  $\g_i$.  This shows that $\g_i$ and $\g_{nra_k+i}$ have the same weights, thus, are isomorphic.

Next assume $\g_i$ and $\g_j$ are isomorphic as $\gk$ modules. 
By Lemma~\ref{afdual}, we may assume $i$ and $j$ are both positive or both negative so say $i$, $j>0$.  A highest root associated to $\g_j$ has the form $\beta=j\alpha_k+\sum_{m\neq k}c_m\alpha_m$ and a highest root associated to $\g_i$ has the form $\gamma=i\alpha_k+\sum_{m\neq k}c'_m\alpha_m$.  Since $\beta$ and $\gamma$ belong to $Q_+$, and since $\gamma-\beta$ restricted to $\hk$ is identically zero, it follows that $\gamma-\beta=n\delta$ for some $n\in{\mathbb Z}$.  In the nontwisted algebras, this is enough to prove the result.

Now suppose $\g$ is type $X_N^{(r)}$, $r=2$ or $3$.  In all cases, the lowest root associated to $\g_1$ is $\gamma_1=\alpha_k$. Comparing Tables Aff~2 and 3 in \cite{kac} to Figure~1 above, we see that $\gamma_2$ is given by Case~1, 6, 7, 8, or 9.  In no case is $\gamma_2=\delta+\gamma_1$.  This finishes the proof for the cases where $r=2$ and $a_k=1$ as $\g_1\not\cong\g_2$ and for $i>2$, $\g_i$ is isomorphic to $\g_1$ or $\g_2$.  When $\g=D_4^{(3)}$, $\gamma_2$ is as in Case~7 of Figure~1.  Using the numbering of roots in \cite{kac}, we find $$\gamma_3=2\alpha_0+5\alpha_1+3\alpha_2=3\delta-(\alpha_0+\alpha_1)=2\delta+(\alpha_1+\alpha_2).$$  (Note that $2\delta+\alpha_k$ is not a root as $\alpha_k$ is long.) Here again, $\gamma_i+\delta$ is not a lowest root for $i=1,2$ or $3$ so the result holds.  Finally, when $\g=E_6^{(2)}$ and $a_k=2$, $\gamma_3$ is  as in Case~5 of Figure~2.   $\gamma_3-\delta$ is not a root (again, $\gamma_3$ is long) and $\gamma_3-2\delta$ is a negative root so no root of the form $\gamma_3+n\delta$ is a lowest positive root for $\g_i$.  A lowest root for $\g_4$ has the form $2\delta-\theta$, $\theta$ a highest root of $\gk$.  Since $\gk=A_3\oplus A_1$, this confirms that no two of $\gamma_1$, $\gamma_2$, $\gamma_3$, or $\gamma_4$ differ by a multiple of $\delta$. The result follows since $\g_i$ is determined by $\g_i\cong \g_{i+4}$.
  $\Box$  

In any affine (or finite type) Lie algebra, $a_k\leq 6$.   Lemmas~\ref{afdual}, \ref{grdel}, and \ref{gigjcon} together thus imply that $\g_1$, $\g_2$, and $\g_3$ determine all other graded pieces of $\g$, either directly or by duality.  

The next theorem gathers some of the results we established in the lemmas.

\begin{thm}\label{afgrd}
\begin{enumerate}
\item The highest weight of $\g_{-1}$ is $\lambda_{-1}=\F{1}{a_k}\delta-\alpha_k$. The highest weight of $\g_1$ is $\F{1}{a_k}\delta-\sigma(\alpha_k)$, where $\alpha_k$ and $\sigma(\alpha_k)$ are dual roots, that is, conjugate under an automorphism of the Dynkin diagram.
\item The highest weight of $\g_{-1}$ is also given by 
$\lambda_{-1}=\sum_{i\neq k} s_i\Lambda_i$ where $s_i$ is the number of edges in $\D$ connecting $\alpha_i$ and $\alpha_k$. $\g_1$ has highest weight $\lambda_{1}=\sum_{i\neq k} s'_i\Lambda_i$ where $s'_i$ is the number of edges in $\D$ connecting $\alpha_i$ and $\sigma(\alpha_k)$.  In particular, $\g_1$ is a generalized cominuscule representation of $\gk$.
\item The grading on an affine Lie algebra $X_N^{(r)}$ determined by a long simple root $\alpha_k$ has period $ra_k$.
\item The grading defined by a long simple root $\alpha_k$ on an affine Lie algebra is completely determined by $\gk$, $\g_1$, $\g_2$ and $\g_3$. 
\end{enumerate} 
\end{thm}

{\bf Proof}  It is clear that $\lambda_{-1}=\pi(-\alpha_k)=\F{1}{a_k}\delta-\alpha_k$ is the highest weight of $\g_1$.  The second statement follows from the fact that $\g_1$ and $\g_{-1}$ are dual representations.  

If $\alpha_i$ is a simple root of $\gk$, $<\lambda_{-1},\alphahat_i>=<-\alpha_k,\alphahat_i>$ is the number of edges between $\alpha_k$ and $\alpha_i$ in $\D$, since $\alpha_k$ is long.  The form of $\lambda_1$ follows by duality. Since there is one component of $\Dk$ for each simple root adjacent to $\alpha_k$ in $\D$, $\g_1$ is generalized cominuscule. 

The rest of the results have already been established.
$\Box$

\subsection{The Finite Case}  Let $\g$ be finite type so that $\h=\hk\oplus\; {\mathbb C} d_k$ and $\h^*=\hkd\oplus\; {\mathbb C}\Lambda_k$ with $d_k=\nu^{-1}(\Lambda_k)$.
It is clear that with respect to the $\alpha_k$-grading on $\g$,
$\g_0=\gk\oplus\; {\mathbb C}d_k.$ The following is an analog to Lemma~\ref{affinefactslem}.

\begin{lem}\label{finitefactslem}
$$\begin{array}{llllr}
(\alphahat_i,d_k)=\delta_{ik} & &(\Lambda_k,\alpha_k)=\F{\ahat_k}{a_k}=1&&\nu(\alphahat_k)=\alpha_k\\ \\

\nu(d_k)=\Lambda_k&&<\alpha_i,d_k>=\delta_{ik}&& \Box \\ \\
\end{array}$$
\end{lem}

Project $\h^*$ onto $\hkd$ by
$\pi:\h^*\rightarrow \hkd$
where $\pi(\alpha_i)=\alpha_i$, for $i\neq k$, and $\pi(\Lambda_k)=0$.  If
$\Lambda_k=\sum_{j}t_j\alpha_j$ we have
$\pi(\alpha_k)=-\sum_{i\neq k}\F{t_i}{t_k}\alpha_i=-\F{1}{t_k}\Lambda_k+\alpha_k.$

\begin{lem}
$\g_1$ is an irreducible $\gk$ module.
\end{lem}

{\bf Proof}  The first part of the proof of Lemma~\ref{afg1ir} applies here.  We do not need a separate argument for the case $a_k=1$ because in the current setting, all roots of $\g$ are real. $\Box$ 

The proof of Lemma~\ref{afdual} applies nearly unchanged to give us the following.

\begin{lem}
$\g_i$ and $\g_{-i}$ are dual representations of $\gk$. $\Box$
\end{lem}

Let $\theta$ be the highest root of $\g$ so that $\theta=\sum_{i>0} a_i\alpha_i$. Let $\Lambda=\pi(\theta)$ so $$\Lambda=\sum_{i\neq k}\left(a_i-\F{t_i}{t_k}\right)\alpha_i.$$

\begin{prop} As $\gk$ modules, $\g_{a_k}\cong V(\Lambda)$. \end{prop} 

{\bf Proof}  $\g_{a_k}$ is a sum of root spaces, among them $\g_\theta$.  Since $\theta$ is the highest root associated to $\g$, it is the highest root associated to $\g_{a_k}$.  The highest weight associated to $\g_{a_k}$ as a $\gk$ module is thus $\Lambda=\pi(\theta)$.  All roots here are real so $\dim\g_{\theta}=1$.  If $v^+$ is a highest weight vector associated to another root space in $\g_{a_k}$, $\g_\beta$, then $\beta$ has the property that $\beta+\alpha_i$ is not a root for any $i\neq k$.  Since $a_k$ is the coefficient of $\alpha_k$ in $\theta$ as well as in $\beta$, we cannot add $\alpha_k$ to $\beta$ and get a root of $\g$ either.  Thus, $\beta+\alpha_i$ is not a root for any simple root $\alpha_i$ associated to $\g$, which is impossible unless $\beta=\theta$, thus, $\g_{a_k}$ is irreducible as a $\gk$ module.  $\Box$

\begin{rem}  If $a_k>1$,  $\dim\g_{a_k-1}=\F{1}{a_k-1}\dim\g_1$. The reader can verify this curious fact case by case for $\g$ finite type but we have no insight beyond the observation itself. There appears to be no analog for affine algebras.
  \end{rem}

\begin{thm}\begin{enumerate}
\item The highest weight of $\g_{-1}$ is $\lambda_{-1}=\F{1}{t_k}\Lambda_k-\alpha_k$, where $\Lambda_k=\sum_it_i\alpha_i$.  The highest weight of $\g_1$ is $\lambda_1=-\F{1}{t'_k}\Lambda'_k+\sigma(\alpha_k)$ where $\Lambda'_k$ is fundamental on $\sigma(\alpha_k)$, the root dual to $\alpha_k$, and $\Lambda'_k=\sum_i t'_i\alpha_i$.  
\item  The highest weight of $\g_{-1}$ is also given by $\lambda_{-1}=\sum_{i\neq k} s_i\Lambda_i$ where $s_i$ is the number of edges in $\D$ connecting $\alpha_i$ and $\alpha_k$.
$\g_1$ has highest weight $\lambda_{1}=\sum_{i\neq k} s'_i\Lambda_i$ where $s'_i$ is the number of edges in $\D$ connecting $\alpha_i$ and $\sigma(\alpha_k)$. In particular, $\g_1$ is generalized cominuscule.  \end{enumerate}
 \end{thm}

{\bf Proof} It is clear that the highest weight associated to $\g_{-1}$ is $\pi(-\alpha_k)=\F{1}{t_k}\Lambda_k-\alpha_k.$  For $i\neq k$, $<\lambda_{-1}, \alphahat_i>=<-\alpha_k,\alphahat_i>$, which is the number of edges between $\alpha_k$ and $\alpha_i$ in $\D$, since $\alpha_k$ is long.  The forms of $\lambda_1$ follow by duality.  As in the affine case, $\gk$ has one simple component per simple root adjacent to $\alpha_k$ in $\D$.  This, along with the form of $\lambda_1$ given in the second statement, gives us that $\g_1$ is generalized cominuscule.
$\Box$  

\subsection{Representations}
We revisit Lemma~\ref{subrootlemma} as a statement about weights.

\begin{prop} \label{subwtprop}
Let $\lambda_i$ be a highest weight of $\g_{i}$ as a $\gk$-module.
If $a_k>1$, then 
$$\lambda_2=2\lambda_1-(\beta_1+\beta_2+\beta_3)$$
for positive roots $\beta_i$ with properties 1-5 in Lemma~\ref{subrootlemma}.
If $a_k\geq t>2$,
$$\lambda_t=t\lambda_{1}-\sum_{i=1}^{2t-1}\beta_i,$$
for positive roots $\beta_i$ arising as in Lemma~\ref{subrootlemma}.  In particular, 
$\lambda_t=\lambda_{t-1}+\lambda_1-\beta_1-\beta_2,$
where $\beta_1, \beta_2\in\;\Deltak_+$ satisfy properties 1-4 in Lemma~\ref{subrootlemma}, {\em mutatis mutandis}.
\end{prop}

{\bf Proof}  If $\gamma_t$ is a lowest root associated to $\g_t$, then $\pi(\gamma_t)$ is a lowest weight of $\g_t$.  Let $w$ be the longest element in the Weyl group of $\gk$:  then $w(\pi(\gamma_t))=\lambda_t$, the highest weight of $\g_t$.  

Lemma~\ref{subrootlemma} gives us $\gamma_2=2\gamma_1+\beta_1+\beta_2+\beta_3$.  Then 
$$w(\pi(\gamma_2))=\lambda_2=2w(\pi(\gamma_1+\beta_1+\beta_2+\beta_3))=2\lambda_1-(\beta_1'+\beta_2'+\beta_3')$$
where the $\beta_i'$s are positive roots.  The longest element in the Weyl group of $\gk$  respects properties 1-5 of Lemma~\ref{subrootlemma}.  Thus, the $\beta_i'$s enjoy those properties as well. The argument for $t>2$ is similar.  
$\Box$

We have established that $\g_1$ and $\g_{-1}$ are irreducible $\gk$ modules.  Next we show that all $\g_i$ are irreducible when $i\not\equiv 0$ mod $ra_k$, in the affine case, $i\neq 0$ in the finite case.  First we need a technical lemma.

If $\alpha$ is a subroot of $\beta$, i.e., $\beta-\alpha\in Q_+$, we write $\alpha\subseteq\beta$. When simple $\alpha_i$ is in the support of $\alpha$, we write $\alpha_i\in\alpha$.

\begin{lem}\label{dilgenlemma}
If $\alpha$ is a positive root of $\g$ and $\alpha_r\subseteq \alpha$, we can write  $\alpha=\alpha_r+\sum_{j=1}^m\alpha_{i_j}$ for $\alpha_{i_j}\;\in\;\Pi$ where, for all $t\;\in\;\{1,\ldots,m\}$, $\alpha_r+\sum_{j=1}^t\alpha_{i_j}\;\in\;\Delta_+.$
\end{lem}

{\bf Proof}  The case where $\alpha_r$ is simple is done in \cite{dil}, Lemma~2.1, so here we assume that the height of $\alpha_r$ is greater than one.

If there is a simple root $\alpha_i\not\in\alpha_r$ with $\beta=\alpha-\alpha_i\;\in\;\Delta_+$, we get the result by induction on the height of $\alpha$ because $\alpha_r\subseteq \beta$.

If no such simple root exists, we proceed by induction on $\mbox{ ht }(\alpha-\alpha_r)$, the case $\mbox{ ht }(\alpha-\alpha_r)=1$ being clear.  If $\hgt(\alpha-\alpha_r)>1$, the result follows once we produce a simple root $\alpha_j\in\alpha-\alpha_r$ with $\alpha_r+\alpha_j\in\Delta_+$. 

If there is a simple root $\alpha_j\in\alpha$, $\alpha_j\not\in\alpha_r$, we claim there must be one that satisfies $\alpha_r+\alpha_j\in\Delta_+$. The support of any root must form a connected component of the Dynkin diagram, which implies there must be $\alpha_j$ that shares an edge of $\D$ with some $\alpha_i\in\alpha_r$.  Then $<\alpha_r,\alphahat_i>\;<0$ implies $\alpha_r+\alpha_j\in\Delta_+$ as claimed.  That leaves us with the case where $\supp \alpha=\supp \alpha_r$. 

If $(\alpha_r,\alpha)<0$, there must be a simple $\alpha_i\in\alpha$ with $\alpha_r+\alpha_i\in\;\Delta_+$.  If  $(\alpha_r,\alpha)=0$, invoke the fact that $\supp \alpha_k=\supp\alpha$.  This, along with the fact that $\hgt \alpha_r>1$ implies there is $\alpha_i\in\supp\alpha_r$ with $\alpha_r-\alpha_i\in\;\Delta_+$, means that there must be $\alpha_j\in\supp\alpha_r=\supp\alpha$ with $\alpha_r+\alpha_j\in\;\Delta_+$.  Finally, suppose $(\alpha_r,\alpha)>0$.  This gives us $\alpha-\alpha_r\in\;\Delta_+$.  Now invoke $\alpha_r+(\alpha-\alpha_r)\in\;\Delta_+$ to get some $\alpha_i\in\supp (\alpha-\alpha_r)$ with $\alpha_r+\alpha_i$ a root. $\Box$ 

\begin{thm}
If $\g=X_N^{(r)}$ is affine and $i\not\equiv 0$ mod $ra_k$, $\g_{i}$ is an irreducible $\gk$ module. If $\g$ is finite type and $i\neq 0$, $\g_i$ is an irreducible $\gk$ module.
\end{thm}

{\bf Proof} To show $\g_i$ is irreducible, we show it has a unique lowest weight.

We show below that the lowest root for $\g_i$ is well-defined.  For now, assume this is so and let $\gamma_i$ be the lowest root in $\g_i$.  If there is a second root $\gamma'_i$ such that $f_j.u=0$ for some nonzero $u$ in $\g_{\gamma'_i}$, and  all $\alpha_j\in\Pik$, note that $\gamma_i\subseteq\gamma'_i$.  Lemma~\ref{dilgenlemma} allows us to write 
$\gamma'_i=\gamma_i+\sum_{j=1}^m\alpha_{i_j}$ where $\alpha_{i_j}\in\;\Pik$ and for all $t\in\{1,\ldots,m\}$, $\gamma_i+\sum_{j=1}^t\alpha_{i_j}\in\;\Delta_+$.  Applying Proposition~3.6 in \cite{kac},  we get $f_{i_1}.\g_{\gamma'_i}\neq 0$.  Since $\gamma_i$ and $\gamma'_i$ have the same coefficient on $\alpha_k$, $\alpha_{i_1}$ must belong to $\Pik$.  Our assumption that $f_{i_1}.u=0$ forces us to conclude that $\dim\g_{\gamma'_i}>1$, thus, that $\gamma'_i$ is imaginary. An argument similar to the one we used for Lemma~\ref{gigjcon} applies to show that there must be some $f_{i_j}$ with $\alpha_{i_j}\in\Pik$ and $f_{i_j}.u\neq 0$.  The contradiction proves that $\g_i$ is irreducible.

Now we show that there is a unique lowest root $\gamma_i$ associated to $\g_i$.

There is only one choice of $\gamma_1=\alpha_k\in\h^*$ so the result holds when $i=1$ thus when $i=-1$.  By duality and the fact that $a_k\leq 6$, it suffices to show that if $i=2$ or $3$ and $ra_k>2$ or $3$ respectively, then for a fixed $\g$ and fixed $\alpha_k$, there is only one choice of $\gamma_i$ from among the root diagrams in Figures~1 and 2 respectively, at least in the affine case.   We deal with those details and leave the rest of the finite type cases to the reader.

A choice of $\gamma_i$ corresponds to a choice of injective mapping of ${\D}(\gamma_i)$, up to automorphism of ${\D}(\gamma_i)$, into $\D$ with its $\alpha_k$ node marked.  We do our accounting by considering Tables Aff~1, 2, and 3 in \cite{kac} that have nodes associated to long roots with labels such that $ra_k>2$.  We compare those to the diagrams from Figures~1 and 2 that inject into ${\D}$ and see that in each case, as long as $i<ra_k$, there is a unique injection of ${\D}(\gamma_i)$ into $\D$, up to automorphism of ${\D}(\gamma_i)$.  

Consider the case $i=2$ where $ra_k>2$.  
$\g$ can be type $E$, $F$, or $D_4^{(3)}$.  Since $D_4^{(3)}$, $E_6^{(2)}$, and $F$ type algebras are not simply laced, it is clear that there is only one way to inject ${\D}(\gamma_2)$ in $\D$ in these cases.  (These correspond respectively to Cases~7, 6, and 1 in Figure~1.)  

Next suppose ${\D}$ is nontwisted type $E$ with the branch node marked.  This is Case~2 in Figure~1.  There are three different ways to inject that diagram into ${\D}$ with the branch node distinguished but these different injections arise from diagram automorphisms of ${\D}(\gamma_2)$.  

If $\alpha_k$ corresponds to the terminal node at the end of the shortest branch of $\D$, ${\D}(\gamma_2)$ is Case~3 in Figure~1.  Ignoring diagram automorphisms, there is only one way to inject this into $\D$.
If $\alpha_k$ is the terminal node at the end of a long branch, then $\g=E_8^{(1)}$ and the diagram from Figure~1 is Case~5.  There is a unique injection of ${\D}(\gamma_2)$ into $\D$.   

If $\g$ is nontwisted $E$ type with some other distinguished node, the relevant diagram in Figure~1 is Case~2.  If ${\D}(\gamma_2)$ is a $D_5$ diagram,  it fits two different ways into $E_7^{(1)}$ and $E_8^{(1)}$ diagrams.  But with a non branch node distinguished on each diagram, there is only one way to inject ${\D}(\gamma_2)$ into ${\D}$.  

This is a complete proof that $\g_2$ is irreducible in case $ra_k>2$.  Next  consider $\g_3$ when $ra_k>3$.  

Here $\D$ can be type $E_6^{(2)}$,  $E_7^{(1)}$ or $E_8^{(1)}$.  Again, since $E_6^{(2)}$ is not simply laced, the uniqueness of the diagram injection is immediate.  (${\D}(\gamma_3)$ is Case~5 in Figure~2.)  In the $E_7^{(1)}$ case the distinguished node is the branch node and ${\D}(\gamma_3)$  comes from Case~1 in Figure~2.  Note that ${\D}(\gamma_3)$ is a type $E_6$ diagram.  There are two ways to inject such a diagram into one of type $E_7^{(1)}$ but the different ways arise from an automorphism of the $E_6$ diagram.  

When the distinguished node on $\D$ is not the branch node, it is one or two nodes away from the branch node.  If one away, it is adjacent to a terminal node or not.  If $\alpha_k$ is between a terminal node and a branch node, ${\D}(\gamma_3)$ comes from Case~2 in Figure~2 so is of type $E_7$.  There is a unique way to inject an $E_7$ diagram into an $E_8^{(1)}$.  If $\alpha_k$ is adjacent to the branch, not adjacent to a terminal node, ${\D}(\gamma_3)$ is again from Case~1 Figure~2.  This also injects uniquely into $\D$.   $\Box$  

\subsection{The Invariant Form and The Casimir Operator}
How is $(.,.)|_{\gk}$ related to the Killing form on $\gk$? 

If $V(\lambda)$ is a generalized cominuscule representation of a semi-simple Lie algebra $\g$ with $\lambda=\sum_i n_i\Lambda_i$ as in Definition~\ref{com}, the {\em degree} of the representation on $\mcl_i$ is $n_i$. 

Let $\mcl_i$ be a simple component of $\gk$.  Let $n_i$ be the degree of $\g_1$ as a generalized cominuscule representation of $\mcl_i$.  If $B(.,.)$ is the standard invariant form on $\gk$ normalized so that the square length of a long root is $2$, then $$(.,.)|_{\mcl_i}=\F{1}{n_i}B(.,.).$$ This maintains relative root lengths in $\g$ on restriction to $\gk$. Conversely, if we start with semi-simple $\gk$ in a generalized cominuscule representation, and construct $\g$ using one of the algorithms, $(.,.)|_{\mcl_i}:=\F{1}{n_i}B(.,.)$ lifts to the invariant form on $\g$ that measures the square length of a long root as $2$.  Note:  (1) long roots of $\gk$ belonging to different simple components may have different lengths; (2) insisting that the square length of long roots of $\g$ is $2$ is not the standard normalization on twisted affine Lie algebras; (3) it remains true that $<\alpha_i,\alphahat_j>=\F{2(\alpha_i,\alpha_j)}{\|\alpha_j\|^2}$ is the number of edges shared by $\alpha_i,\alpha_j\in\Pik$ and when $\alpha_i$ is long.  

Let $\lambda_j$ be a highest weight of $\g_j$. 

\begin{lem} \label{lam1}\begin{enumerate} \item In the affine case, $(\lambda_1,\lambda_1)=2.$  In the finite case, $(\lambda_1,\lambda_1)=2-\F{1}{t_k}$, where $\Lambda_k=\sum_it_i\alpha_i$.
\item For $\alpha_i\in \Pik$, $(\lambda_{1},\alpha_i)=1$ or $0$. \end{enumerate}\end{lem}

{\bf Proof}  \begin{enumerate}\item In the affine case, 
$(\lambda_1,\lambda_1)=\left\|\F{1}{a_k}\delta-\alpha_k\right\|^2=\|\alpha_k\|^2=2.$
The proof in the finite case is similar.
\item Suppose  $<\lambda_{1},\alphahat_i>=n$ so that $\|\alpha_i\|^2=2/n$.  Then 
$(\lambda_{1},\alpha_i)=<\lambda_{1},\alphahat_i>\F{\|\alpha_i\|^2}{2}=1.$ $\Box$ \end{enumerate}

If $V(\lambda)$ is an irreducible highest weight representation of a semi-simple Lie algebra $\g$, the Casimir operator $C$ takes the scalar value $c=(\lambda,\lambda)+2(\rho,\lambda)$, where $\rho$ is the sum of fundamental weights associated to $\g$.

In what follows, $c_i$ is the scalar value of the Casimir operator on $\g_i$ when $\g_i$ is irreducible and $\rho$ is the sum of the fundamental weights associated to $\gk$.

\begin{prop}\label{casprop} $\cT_1$ is a Casimir eigenspace.  In the affine case, $C|_{\cT_1}=2c_1$.  In the finite case, $C|_{\cT_1}=2c_1-\F{2}{t_k}.$\end{prop}

{\bf Proof} A highest weight of $\cT_1$ has the form $2\lambda_1-\alpha_i$ where $\alpha_i$ is a simple root of $\gk$ with $<\lambda_1,\alphahat_i>\neq 0$.  On the associated irreducible component of $\cT_1$, we have
$$C=(2\lambda_1-\alpha_i,2\lambda_1-\alpha_i)+2(\rho,2\lambda_1-\alpha_i)=2c_1+2(\lambda_1,\lambda_1)-4(\lambda_1,\alpha_i)+\|\alpha_i\|^2-2(\rho,\alpha_i).$$
In the affine case,  $(\lambda_1,\lambda_1)=2$ by Lemma~\ref{lam1}.  We also have 
$$(\rho,\alpha_i)=<\rho,\alphahat_i>\F{\|\alpha_i\|^2}{2}=\F{\|\alpha_i\|^2}{2}.$$
Since $(\lambda_1,\alpha_i)=1$, we get
$$C=2c_1+4-4=2c_1.$$

When we do the calculation for the finite case, the only change is $\|\lambda_1\|^2=2-\F{1}{t_k}$.  $\Box$

Let $\beta_i$s be as in Lemma~\ref{subrootlemma} and let $n$ be the degree of the representation $\g_2$ on a given component of $\gk$.  Perusal of the diagrams in Figure~1 reveals two important data:  (1) the $\beta_i$s are all long roots in their respective components of $\Dk$; and (2) in Cases7-9, those for which $n>1$, the $\beta_i$s  belong to a single component of $\Dk$.  In particular, $\|\beta_i\|^2=2/n$, where $n$ is well-defined in each case. 

\begin{cor} 
\begin{enumerate}\item Suppose $ra_k>2$.  In the affine case, the Casimir operator acts on $\g_2$ as the scalar
$c_2=2c_1-8+\F{6}{r}-\F{2}{r}\sum_{i=1}^3\hgt\beta_i.$ 
\item Suppose $ra_k\geq 2$.  
In the finite case, the Casimir operator acts on $\g_2$ as the scalar
$c_2=2c_1-8-\F{2}{t_k}+\F{6}{n}-\F{2}{n}\sum_{i=1}^3 \hgt\beta_i.$\end{enumerate}
\end{cor}
{\bf Proof} Consider that
$$c_2=(2\lambda_1-(\beta_1+\beta_2+\beta_3),2\lambda_1-(\beta_1+\beta_2+\beta_3))+2(\rho, 2\lambda_1-(\beta_1+\beta_2+\beta_3))$$$$=4(\lambda_1,\lambda_1)+4(\rho,\lambda_1)-4(\lambda_1,\beta_1+\beta_2+\beta_3)+\|\beta_1+\beta_2+\beta_3\|^2-2(\rho,\beta_1+\beta_2+\beta_3)$$$$=2c_1+2(\lambda_1,\lambda_1)-4\cdot 3+\|\beta_1+\beta_2+\beta_3\|^2-2(\rho,\beta_1+\beta_2+\beta_3).$$
We note
$$(\rho,\beta_i)=<\rho,\betahat_i>\F{\|\beta_i\|^2}{2}=\F{\hgt \beta_i}{n}$$ giving us
$$c_2=2c_1+2(\lambda_1,\lambda_1)-12+\F{6}{n}-\F{2}{n}\sum_{i=1}^3 \hgt\beta_i.$$
Suppose $ra_k>2$ in the affine case.  Perusal of the tables in \cite{kac} reveals that $n=r$ in these cases.  By Lemma~\ref{lam1} we have
$$c_2=2c_1-8+\F{6}{r}-\F{2}{r}\sum_{i=1}^3\hgt\beta_i.$$
If $a_k\geq 2$ in the finite case, we apply Lemma~\ref{lam1} to get
$$c_2=2c_1+4-\F{2}{t_k}-12+\F{6}{n}-\F{2}{n}\sum_{i=1}^3 \hgt\beta_i$$
$$=2c_1-8-\F{2}{t_k}+\F{6}{n}-\F{2}{n}\sum_{i=1}^3 \hgt\beta_i. \Box $$

Let $Y_{1,t}$ be the Cartan product of $\g_1$ and $\g_t$, that is, $Y_{1,t}$ is irreducible with highest weight $\lambda_1+\lambda_t$.  A corollary of Proposition~\ref{subwtprop} follows.

\begin{cor}\label{lam1lamt}
\begin{enumerate}
\item If $\g$ is affine, $(\lambda_1,\lambda_t)=1$, for $t\in\{2,\ldots,a_k\}$.
If $\g$ is finite type, $(\lambda_1,\lambda_t)=1-\F{t}{t_k}$.
\item If $\g$ is affine, the value of the Casimir operator on $Y_{1,t}$ is $c_1+c_t+2$.  If $\g$ is finite type, the value of the Casimir on $Y_{1,t}$ is $c_1+c_t+2-\F{2t}{t_k}$. 
\item If $\g$ is affine, the value of the Casimir operator on a component of $\cT_t$ with highest weight $\lambda_1+\lambda_t-\alpha$ is $c_1+c_t+2-2(\rho,\alpha)$. If $\g$ is finite type, the Casimir on such a component is $c_1+c_t+2-\F{2t}{t_k}-2(\rho,\alpha)$. 
\end{enumerate}
\end{cor}

{\bf Proof} \begin{enumerate}\item Proposition~\ref{subwtprop} along with Lemma~\ref{lam1} gives us $$(\lambda_1,\lambda_t)=t(\lambda_1,\lambda_1)-(2t-1)=2t-2t+1=1$$ when $\g$ is affine and when $\g$ is finite type,
we adjust according to $(\lambda_1,\lambda_1)=2-\F{1}{t_k}$.
\item The Casimir operator takes the following value on $Y_{1,t}$: $$(\lambda_1+\lambda_t,\lambda_1+\lambda_t)+2(\lambda_1+\lambda_t,\rho)=c_1+c_t+2(\lambda_1,\lambda_t).$$ Apply statement (1) now to get statement (2).
\item If $\g$ is affine type, the Casimir on a given component of $\cT_t$ is
$$c_1+c_t+2-2(\alpha,\lambda_1+\lambda_t)+(\alpha,\alpha)-2(\rho,\alpha).$$  We have $(\lambda_1,\alpha)=(\lambda_t,\alpha)=1$ so the Casimir is $c_1+c_t-2-2(\rho,\alpha)$.  In the finite case, adjust the calculation to reflect the result in part (2). $\Box$ 
\end{enumerate}

$\g_2$ is a distinguished submodule of $\Lambda^2\g_1$ and for $i>2$, 
$\g_i$ is a distinguished submodule of $\g_{1}\otimes \g_{i-1}$.  The bracket is then a projection from $\g_{1}\otimes \g_{i-1}$ onto $\g_i$. Proposition~\ref{subwtprop} tells us what the kernel of the projection is.

\begin{cor} The bracket as defined on $\Lambda^2\g_1$ is identically zero on $\cT_1$.  
If $2\leq i \leq a_k-1$, the bracket defined on $\g_1\otimes \g_i$ is identically zero on $\cT_i$. $\Box$\end{cor}

\section{Proof of the Algorithm, Part I}
It is convenient at this point to treat the Affine Algorithm and the Finite Algorithm as one process.   Our goals in this section are: (1) to show that the algorithms are effective, either aborting or terminating successfully with $j\leq 6$; and (2) to establish that there is a one-to-one correspondence between the graded $\gk$ modules produced by the algorithms and lsn-graded affine/finite type Lie algebras. (Note that $A^{(1)}_n$ and $C_n^{(1)}$ are the only affine or finite type Lie algebras without lsn-gradings.)   In the next section, we establish that if $\g$ is produced by one of the algorithms, then it has the expected Lie algebra structure.

The effect of a successful run of either algorithm is to append a long node to $\Dk$.  When $V$ is degree $n$ on a component of $\gk$, the associated node of $\Dk$ shares $n$ edges with the appended node. We use this idea to catalog the graded modules produced by the algorithm.

When considering different input modules, we make no distinction between duals or between different modules that yield dual outputs, for example, the two half-spin representations of $D_n$.
  
\begin{lem}\label{d41} The algorithm aborts with $j=1$ in cases where $\gk$ has more than four components. The only case in which it tolerates four components is that for which $V=U_1\otimes\ldots\otimes U_4$, where $U_i$ is the standard (two dimensional) representation of $A_1$.  This case terminates with $j=2$ and corresponds to $\alpha_3$-graded $D_4^{(1)}$.\end{lem}

{\bf Proof} Suppose $\gk=\mcl_1\oplus\ldots\oplus\mcl_5$ and $V=U_1\otimes \ldots\otimes U_5$, $U_i$ generalized cominuscule on $\mcl_i$.  $\Lambda^2 V$ has sixteen (not necessarily irreducible) components, each a product of one, three, or five $\Lambda^2 U_i$s, with, respectively, four, two, or zero $S^2 U_j$s.  A weight argument confirms that $\cT_1$ is contained in the sum of the five components that have one $\Lambda^2U_i$ and four $S^2 U_j$ factors.  This leaves at least eleven more irreducible components, too many to comprise either an irreducible $\cT_1^c$ or the adjoint representation of $\gk$.
If $\gk$ has more than five components, the number of components of $\Lambda^2V$ is even farther beyond the maximum tolerated by the algorithm.  So in all cases where $\gk$ has five or more simple components, the algorithm aborts with $j=1$.

Next suppose $\gk$ has four components.  This time, $\Lambda^2 V$ has eight not necessarily irreducible components:  four each consisting of the tensor product of one $\Lambda^2 U_i$ with three $S^2 U_j$s (this sum contains $\cT_1$), and four each consisting of the tensor product of three $\Lambda^2 U_i$s with one $S^2 U_j$ (these are contained in $\cT_1^c$.) Since the sum of four submodules contained in $\cT_1^c$ is not irreducible, the algorithm aborts with $j=1$ unless $\cT_1^c$ coincides with the adjoint representation of $\gk$.  Note in particular that the algorithm aborts if any of the eight components of $\Lambda^2 V$ is not irreducible.  If the algorithm does not abort, the submodule of $\Lambda^2 V$ consisting of products of one symmetric square and three exterior squares must be of the form $\adl \otimes {\mathbb C}\otimes {\mathbb C}\otimes {\mathbb C}$, up to ordering of the factors.  This establishes that $\Lambda^2 U_i={\mathbb C}$, thus, that $U_i$ is the two dimensional representation of ${\it sl}_2({\mathbb C})$.  The algorithm terminates successfully then with $j=2$. $\Box$ 

\subsection{$V=U_1\otimes U_2\otimes U_3$}

Suppose $V=U_1\otimes U_2\otimes U_3$, where $U_i$ is a generalized cominuscule representation of ${\mfl}_i$.  For the remainder of this section, we use $n_i$ to designate the rank of $\mfl_i$.  To make accounting easier, assume $n_1\leq n_2\leq n_3$.  We have $$\Lambda^2 V=\Lambda^2U_1\otimes S^2U_2\otimes S^2U_3\oplus S^2U_1\otimes \Lambda^2U_2\otimes S^2U_3\oplus S^2U_1\otimes S^2U_2\otimes\Lambda^2U_3$$$$\oplus\Lambda^2 U_1\otimes\Lambda^2U_2\otimes\Lambda^2U_3.$$ $\cT_1^c$ always contains $\Lambda^2 U_1\otimes\Lambda^2U_2\otimes\Lambda^2U_3$. If some $S^2 U_i$ or $\Lambda^2 U_i$ is not irreducible, $\cT_1^c$ has other components as well.  $S^2 U_i$ is irreducible only if $U_i$ is the standard representation of $A_n$.  In this case, $\Lambda^2 U_i$ is irreducible, as well.  Thus, if $\cT_1^c$ has components other than $\Lambda^2 U_1\otimes\Lambda^2U_2\otimes\Lambda^2U_3$, it is because $S^2U_i$ is not irreducible.  

 If $S^2U_i$ is not irreducible and the algorithm advances beyond $j=1$, then  $\Lambda^2U_1\otimes \Lambda^2U_2\otimes \Lambda^2U_3\cong \mbox{ adj }{\mfl}_i$, for some $i$. It follows that $\Lambda^2U_j={\mathbb C}$ for $j\neq i$.  This forces ${\mfl}_i=A_1$ for $i=1,2$.  

We use LiE notation \cite{lie} for irreducible modules, that is, we identify an irreducible representation with the coordinate vector of its highest weight, using the basis of fundamental weights.

\begin{lem} Suppose $\gk=A_1\oplus A_1\oplus \mfl_3$ and $V=U_1\otimes U_2\otimes U_3$.
\begin{enumerate}  \item If $U_1$ and $U_2$ are copies of the standard representation of $A_1$, the algorithm tolerates for $U_3$ only one of the following types of representations:
\begin{enumerate}
\item $\mfl_3=A_n$, $U_3=[1,0,\ldots,0]$ or, $\mfl_3=A_3$, $U=[0,1,0]$, or $\mfl_3=A_1$, $U_3=[2]$.  In these cases, respectively, the algorithm terminates with $j=2$, $\g$ corresponding to $D_{n+3}$ with $\alpha_{n+1}$-grading; $D_5^{(1)}$ with $\alpha_2$-grading; $B_3^{(1)}$ with $\alpha_2$-grading;
\item $\mfl_3=D_n$, $U_3=[1,0,\ldots,0]$; the algorithm terminates with $j=2$, $\g$ corresponding to $\alpha_2$-graded $D_{n+2}^{(1)}$;
\item $\mfl_3=B_n$, $U_3=[1,0,\ldots,0]$; the algorithm terminates with $j=2$, $\g$ corresponding to $\alpha_2$-graded $B_{n+2}^{(1)}$.
\end{enumerate} 
\item If $U_1=[2]$, and $U_2=[1]$, the algorithm aborts with $j=2$ unless $U_3$ is the standard representation of $A_1$.  
\item If either $U_1$ or $U_2$ is $[n]$, $n>2$, the algorithm aborts with $j=2$.
\end{enumerate}

\end{lem}

{\bf Proof}  \begin{enumerate} \item Suppose $U_1$ and $U_2$ are copies of the standard representation of $A_1$.  We have $$\Lambda^2 V= {\mathbb C}\otimes S^2 U_2\otimes S^2 U_3\oplus S^2 U_1\otimes {\mathbb C}\otimes S^2 U_3\oplus S^2 U_1\otimes S^2 U_2\otimes \Lambda^2 U_3\oplus {\mathbb C}\otimes {\mathbb C}\otimes \Lambda^2 U_3.$$  Note that $S^2 U_i=\mbox{ adj }A_1$, for $i=1,2$.
 
 If $\Lambda^2 U_3$ is not irreducible, the algorithm aborts with $j=1$:  it only tolerates a reducible $\cT_1^c$ if $\cT_1^c$ is the adjoint representation of $\gk$, in particular, there can be no more than one irreducible component of $\cT_1^c$ per simple summand of $\gk$.  This means $\mfl_3\neq C_n$.

Suppose $\Lambda^2 U_3$ is irreducible.  If $S^2 U_3$ has more than one nontrivial component, the algorithm aborts with $j=1$: in this case, ${\mathbb C}\otimes \mbox{ adj }A_1\otimes S^2 U_3$ contributes to $\cT_1^c$ something other than the adjoint representation of $A_1$, the net effect being $\cT_1^c$ is neither irreducible nor the adjoint representation of $\gk$.  This guarantees that $\mfl_3\neq E_n$. 

We can restrict attention to representations of type $A$, $B$,  and $D$ algebras with $\Lambda^2U_3$ irreducible, and $S^2 U_3$ with no more than one nontrivial irreducible component.

\begin{enumerate}\item Let $\mfl_3=A_n$.  If $U_3=[1,0,\ldots,0]$, it is easy to verify that the algorithm terminates successfully with $j=2$.
Suppose $U_3=[0,1,0,\ldots,0]$.  If $n\geq 3$, $S^2 U_3=[0,2,0,\ldots,0]+[0,0,0,1,0,\ldots,0]$ so when $n>3$, the algorithm aborts with $j=1$.  When $n=3$, the second representation is $\mathbb C$.  In this case, $\Lambda^2 U_3=\mbox{ adj }A_3$.  The algorithm tolerates this case, designates $V_2=A_1\oplus A_1 \oplus A_3$, and thus terminates successfully at $j=2$. 

If $U_3$ is a higher exterior power of the standard representation of $A_n$, $S^2 U_3$ has more than one nontrivial irreducible submodule so the algorithm aborts in these cases with $j=1$.  The same observation applies to higher degree cominuscule representations of $A_n$, except when $U_3=\mbox{ adj }A_1$.  Here, the algorithm yields $V_2=A_1\oplus A_1\oplus A_1$, terminating successfully with $j=2$.

\item Let $\mfl_3=D_n$.  If $U_3=[1,0,\ldots,0]$, $S^2 U_3=[2,0,\ldots,0]+\mbc$ and $\Lambda^2 U_3=\mbox{ adj }D_n$.  The algorithm designates $V_2=A_1\oplus A_1\oplus D_n$, terminating successfully with $j=2$.

Suppose $U_3=\mbox{ adj }D_n$.  In this case, $S^2 U_3$ has more than one nontrivial component so the algorithm aborts with $j=1$.  

If $n=4$, we can take $U_3=[0,0,1,0]$ because then $S^2 U_3= [0,0,2,0]+\mbc$ and $\Lambda^2 U_3=\mbox{ adj }D_4$. If $U_3$ is a half-spin representation when $n>4$, there is a second nontrivial component of $S^2 U_3$ which makes $\cT_1^c$ neither irreducible nor adjoint.  In these cases, the algorithm aborts with $j=1$.

If $U_3$ is a higher degree cominuscule representation of $D_n$, $S^2 U_3$ has more than one nontrivial component and the algorithm aborts with $j=1$.

\item The cases where $U_3$ is a  representation of $B_n$ are similar to those where $U_3$ is a fundamental representation of $D_n$. We leave details to the reader.  \end{enumerate}

\item Now suppose $U_1=[2]$ and  $U_2=[1]$.  In this case $S^2 U_1=[4]+\mbc$ so $$\cT_1^c={\mathbb C}\otimes {\mathbb C}\otimes S^2U_3\oplus {\mathbb C}\otimes \mbox{ adj }A_1\otimes \Lambda^2U_3\oplus \mbox{ adj }A_1\otimes {\mathbb C}\otimes \Lambda^2 U_3.$$  The algorithm aborts unless this is the adjoint representation of $\gk$.  That happens precisely when $\Lambda^2 U_3$ is trivial and $S^2 U_3$ is the adjoint representation of the third simple component of $\gk$.  That is, $U_3$ must be the standard representation of $A_1$. 

\item Suppose $U_1$ and $U_2$ are copies of $\mbox{ adj }A_1$.  In this case, $\cT_1^c$ has at least three irreducible components, one of which is $\mbox{ adj }A_1\otimes \mbox{ adj }A_1\otimes\Lambda^2U_3$.  This forces the algorithm to abort.  Higher symmetric powers of the standard representation of $A_1$ yield similar results. $\Box$
\end{enumerate}

Next we consider what happens when $\cT_1^c=\Lambda^2 U_1\otimes\Lambda^2 U_2\otimes \Lambda^2 U_3$.  As per remarks preceding the theorem, $\mfl{_i}$ is type $A$ for all $i$ and $U_i=[1,0,\ldots,0]$.
In all these cases, the algorithm advances to $j=2$ with $V_2=\Lambda^2 U_1\otimes\Lambda^2 U_2\otimes \Lambda^2 U_3$.  Let $C_{1,i}(U_j)$ be the Cartan product of $U_j$ and $\Lambda^i U_j$.    We have
$$V_1\otimes V_2=(C_{1,2}(U_1)\oplus\Lambda^3U_1)\otimes(C_{1,2}(U_2)\oplus\Lambda^3U_2)\otimes(C_{1,2}(U_3)\oplus\Lambda^3U_3),$$  so that $$\cT_{2}^c=C_{1,2}(U_1)\otimes\Lambda^3 U_2\otimes \Lambda^3U_3\oplus \Lambda^3U_1\otimes C_{1,2}(U_2)\otimes \Lambda^3U_3\oplus$$$$\Lambda^3U_1\otimes \Lambda^3U_2\otimes C_{1,2}(U_3)\oplus\Lambda^3U_1\otimes\Lambda^3U_2\otimes \Lambda^3U_3.$$  The algorithm directs that one of the following happens.
\begin{enumerate}
\item If $\cT_{2}^c=\{0\}$ the algorithm terminates successfully with $j=2$.
\item If $\cT_{2}^c$ is irreducible, the algorithm advances to $j=3$ with $V_3=\cT_{2}^c$.
\item If $\cT_{2}^c\cong\adg\oplus{\mathbb C}$, $V_3=\adg$ and the algorithm terminates successfully with $j=3$.
\end{enumerate}

\begin{lem} $\cT_{2}^c=\{0\}$ if and only if $U_1$ and $U_2$ are standard representations of $A_1$.  In this case, the algorithm terminates successfully with $\g$ corresponding to $\alpha_{n_3+1}$-graded $D_{n_3+3}$.
\end{lem}

{\bf Proof} We work under the assumption that $n_1\leq n_2\leq n_3$ so $\cT_{2}^c=\{0\}$ precisely when $\Lambda^3 U_i=\{0\}$ for $i=1,2$, which happens only if $n_i=1$, for $i=1,2$.  $\Box$ 

Now we work under the assumption that the algorithm has advanced to $j=2$.

\begin{lem} $\cT^c_2$ is irreducible only if $U_1=[1]$ and $n_2>1$.\end{lem}

{\bf Proof}  $\cT^c_2$ is irreducible precisely when $\Lambda^3U_i=\{0\}$ for exactly one $i$, thus when $i=1$.  This means $\mfl_{n_1}=A_1$.  $\Box$ 

\begin{lem} $\cT^c_2=\adg +\mbc$ if and only if $U_1=U_2=U_3$ is the standard representation of $A_2$.  In this case the algorithm terminates successfully with $j=3$, and $\g$ corresponding to  $\alpha_4$-graded $E_6^{(1)}$. \end{lem}

{\bf Proof}  $\cT^c_2=\adg+\mbc$ only if $\Lambda^3U_i={\mathbb C}$ for all $i$, which implies the result.$\Box$ 

We have proved the following.

\begin{prop}  When $V=U_1\otimes U_2\otimes U_3$, the algorithm advances beyond $j=3$ only if each $U_i$ is the standard representation of $A_{n_i}$ with $n_1=1$, and $n_3\geq n_2\geq 2$.  In this case we have $V_3=U_1\otimes\Lambda^3U_2\otimes\Lambda^3U_3$. $\Box$ \end{prop}

Suppose $V_3$ is irreducible.  We have $$V_1\otimes V_3=(\mbox{ adj }A_1\oplus{\mathbb C})\otimes (C_{1,3}(U_2)\oplus\Lambda^4 U_2)\otimes(C_{1,3}(U_3)\oplus\Lambda^4 U_3).$$  This gives us $$\cT_{3}^c=\mbox{ adj }A_1\otimes\Lambda^4 U_2\otimes \Lambda^4U_3\oplus {\mathbb C}\otimes C_{1,3}(U_2)\otimes \Lambda^4U_3\oplus$$$$ {\mathbb C}\otimes \Lambda^4U_2\otimes C_{1,3}(U_3)\oplus{\mathbb C}\otimes\Lambda^4U_2\otimes \Lambda^4U_3.$$  The algorithm directs that one of the following happens.
\begin{enumerate}
\item If $\cT_{3}^c=\{0\}$ the algorithm terminates successfully with $j=3$.
\item If $\cT_{3}^c$ is irreducible, the algorithm advances to $j=4$ with $V_4=\cT_{3}^c$.
\item If $\cT_{2}^c\cong\adg\oplus{\mathbb C}$, $V_4=\adg$ and the algorithm terminates successfully with $j=4$.
\end{enumerate}

\begin{lem}  Suppose the algorithm has advanced to $j=3$. $V_3$ is irreducible if and only if $U_1$ is the standard representation of $A_1$, and each of $U_2$ and $U_3$ is the standard representation of $A_2$.  In this case, the algorithm terminates successfully with $j=3$ and $\g$ corresponding to $\alpha_4$-graded $E_6$.
\end{lem}

{\bf Proof}  $V_4=\{0\}$ precisely when $\Lambda^4 U_2=\Lambda^4U_3=\{0\}$ which happens precisely when $A_{n_i}=A_2$, $i=2,3$.  $\Box$ 

\begin{lem} Suppose the algorithm has advanced to $j=4$.  $V_4$ is irreducible precisely when $U_2$ is the standard representation of $A_2$ and $n_3>2$.  \end{lem}

{\bf Proof}  $V_4$ is irreducible precisely when $\Lambda^4U_2=\{0\}$ and $\Lambda^4U_3\neq\{0\}$.  The result follows from there. $\Box$ 

\begin{lem} If the algorithm advances to $j=4$, $V_4=\adg$ precisely when each of $U_2$, $U_3$ is the standard representation of $A_3$.  In this case, the algorithm terminates successfully to produce $\g$ corresponding to $\alpha_4$-graded $E_7^{(1)}$. \end{lem}

{\bf Proof}  $V_4=\adg$ precisely when ${\mathbb C}\otimes\Lambda^4 U_2\otimes\Lambda^4 U_3\cong {\mathbb C}$, that is, when $\Lambda^4 U_2\cong\Lambda^4U_3={\mathbb C}$, that is, when $U_2$ and $U_3$ are copies of the standard representation of $A_3$. $\Box$ 

Suppose the algorithm has advanced to $j=4$ and that $V_4={\mathbb C}\otimes U_2\otimes\Lambda^4 U_3$  is irreducible, with $U_2$ the standard representation of $A_2$.  Since $V_1\otimes V_4=U_1\otimes(S^2U_2\oplus U_2^*)\otimes(C_{1,4}\oplus\Lambda^5U_3)$,
\begin{equation}\label{8} \cT_4^c=U_1\otimes U_2^*\otimes\Lambda^5U_3.\end{equation}  This implies the following.

\begin{lem}  When $V=U_1\otimes U_2\otimes U_3$, the algorithm cannot produce $V_5=\adg$. $\Box$ \end{lem}

\begin{lem}  Suppose the algorithm advances to $j=4$ with $V_4$ irreducible. $\cT_4^c=\{0\}$ precisely when $U_1$ is the standard representation of $A_1$, $U_2$ the standard representation of $A_2$, and $U_3$ the standard representation of $A_3$.  In this case, the algorithm terminates successfully with $j=3$ and $\g$ corresponding to $\alpha_4$-graded $E_7$.\end{lem}

{\bf Proof}  $\cT_4^c=\{0\}$ precisely when $\Lambda^5 U_3=\{0\}$, which happens if and only if $U_3$ is the standard representation of $A_3$.  $\Box$ 

\begin{lem}  If $V=U_1\otimes U_2\otimes U_3$ and  for $i=2,3,4,5$ $V_i$ are defined and irreducible, then $U_1$ is the standard representation of $A_1$, $U_2$ the standard representation of $A_2$, and $U_3$ the standard representation of $A_n$ for some $n>3$.  In this case, $V_5=U_1\otimes U_2^*\otimes \Lambda^5 U_3.$ \end{lem}

{\bf Proof}  The proof follows Eq.~(\ref{8}). $\Box$ 

Suppose the algorithm has advanced to $j=5$ with $V_5$ irreducible, in particular, with $n_3>3$. We have $$V_1\otimes V_5=(\mbox{ adj }A_1\oplus{\mathbb C})\otimes(\mbox{ adj }A_2\oplus{\mathbb C})\otimes (C_{1,5}(U_3)\oplus \Lambda^6U_3).$$  If $n_3=4$, $C_{1,5}(U_3)=U_3$ and $\Lambda^6U_3=\{0\}$ but in any case, the following is never zero $$\cT_5^c=\mbox{ adj }A_1\otimes{\mathbb C}\otimes\Lambda^6U_3\oplus {\mathbb C}\otimes\mbox{ adj }A_2\otimes \Lambda^6U_3\oplus$$$$ {\mathbb C}\otimes{\mathbb C}\otimes C_{1,5}(U_3)\oplus {\mathbb C}\otimes{\mathbb C}\otimes \Lambda^6U_3.$$

\begin{lem} $V_6$ is irreducible if and only if $U_3$ is the standard representation of $A_4$.  In this case, the algorithm terminates successfully with $j=6$ and $\g$ corresponding to $\alpha_4$-graded $E_8$. $\Box$ \end{lem}

\begin{lem} Suppose the algorithm has advanced to $j=5$ and assume $n_3>4$.  The algorithm aborts with $j=5$ unless $n_3=5$, in which case it terminates successfully with $j=6$ and $\g$ corresponding to $\alpha_4$-graded $E_8^{(1)}$. \end{lem}

{\bf Proof} If $n_3>4$, $\cT_5^c$ has at least two irreducible components so the algorithm aborts with $j=5$ unless $\cT_5^c=\adg\oplus{\mathbb C}$.  This happens, in turn, only if $\Lambda^6U_3={\mathbb C}$, that is, if and only if $U_3$ is the standard representation of $A_5$. $\Box$ 

The following summarizes what the algorithm produces when $\gk$ has three simple components.

\begin{prop}  When $\gk$ has more than two simple components the algorithm either aborts with $j\leq 5$, or it terminates successfully with $j\leq 6$ and $\g$ corresponding to one of the following types:  $\alpha_2$-graded $B_n^{(1)}$, $\alpha_{n-2}$-graded $D_n$, $\alpha_2$-graded $D_n^{(1)}$, $\alpha_4$-graded $E_n$, or $\alpha_4$-graded $E_n^{(1)}$.$\Box$ \end{prop}

\subsection{$V=U_1\otimes U_2$}
$U_i$ is a generalized cominuscule representation of finite type $\mfl_i$, $V=U_1\otimes U_2$, and $\gk=\mfl_1\oplus\mfl_2$.  We have $\Lambda^2V=\Lambda^2 U_1\otimes S^2 U_2\oplus S^2 U_1\otimes \Lambda^2 U_2.$  The case where $\cT_1^c=\{0\}$ is part of the Minuscule Algorithm \cite{lan}.  To complete our accounting we note that $\cT_1^c=\{0\}$ if and only if $S^2 U_i$ is irreducible for $i=1,2$, that is, if and only if $U_1$ and $U_2$ are standard representations of type $A$ algebras. In such a case, the algorithm terminates successfully with $j=1$ and $\g$ corresponding to $\alpha_k$-graded $A_n$, $1<k<n$.

$\cT_1^c$ is irreducible if and only if $\Lambda^2 U_i$ and $S^2 U_1$ are irreducible and $S^2 U_2$ has two irreducible components.  $U_1$ must be a standard representation of a type $A$ algebra since all other generalized cominuscule representations of finite type algebras have a reducible symmetric square. For the duration of this discussion, let $\mfl_1=A_n$ and $U_1=[1,0,\ldots,0]$.   Choosing $U_2$ from the following list guarantees that $\cT_1^c$ is irreducible.  The list is exhaustive.

\begin{enumerate}
\item $\mfl_2=A_m$, $U_2=[2,0,\ldots,0]$;
\item $\mfl_2=A_m$, $U_2=[0,1,0,\ldots,0]$;
\item $\mfl_2=B_m$ or $D_m$, $U_2=[1,0,\ldots,0]$;
\item $\mfl_2=D_5$, $U_2=[0,0,0,0,1]$;
\item $\mfl_2=E_6$, $U_2=[1,0,0,0,0,0]$. \end{enumerate}

The next sequence of lemmas analyzes these cases in order.  When we encounter twisted affine Lie algebras, we number the simple roots as in \cite{kac}.  

 \begin{lem} Let $\mfl_2=A_m$ and $U_2=[2,0,\ldots,0]$.  The algorithm terminates successfully if and only if one of the following holds. \begin{enumerate}
\item $m=1$, in which case $j=2$ and $\g$ corresponds to $\alpha_{n+1}$-graded $B_{n+2}$; 
\item $n=1$ or $2$ and $m=2$, in which case $j=2$ and $\g$ corresponds to $\alpha_2$-graded $F_4$ or $F_4^{(1)}$, respectively; or
\item $n=1$ and $m=3$, in which case $j=4$ and $\g$ corresponds to $\alpha_3$-graded $E_6^{(2)}$.
 \end{enumerate}
Otherwise, the algorithm aborts with $j\leq 3$.
\end{lem}

{\bf Proof} Under the hypotheses we have $S^2 U_2=[4,0,\ldots, 0]\oplus [0,2,0, \ldots, 0]$. If $m=1$, then $[0,2,0, \ldots, 0]=[0]={\mathbb C}$.  Note that $\cT_1^c=\Lambda^2U_1\otimes [0,2,0,\ldots, 0]$.  This is always nonzero and irreducible so the algorithm always advances to $j=2$.  Next consider $V_1\otimes V_2$.

If $m=1$, $V_2=\Lambda^2U_1\otimes{\mathbb C}$ so that $V_1\otimes V_2=(C_{1,2}(U_1)\oplus \Lambda^3U_1)\otimes(U_2\otimes{\mathbb C}).$ The algorithm terminates successfully with $j=2$ as $\cT_2^c=\{0\}$.  

For $m\geq 2$, $V_2\otimes [0,2,0,\ldots,0]=[2,2,0,\ldots,0] \oplus [1,1,1,0,\ldots,0] \oplus[0,0,2,0,\ldots,0].$  This gives us $V_1\otimes V_2=(C_{1,2}(U_1)\oplus \Lambda^3U_1)\otimes([2,2,0,\ldots,0] \oplus [1,1,1,0,\ldots,0] \oplus[0,0,2,0,\ldots,0]),$ so that $$\cT_2^c=C_{1,2}(U_1)\otimes [0,0,2,0,\ldots,0]\oplus\Lambda^3U_1\otimes [1,1,1,0,\ldots,0]$$\begin{equation}\label{9} \oplus\Lambda^3U_1\otimes [0,0,2,0,\ldots,0].\end{equation} 

Expression (\ref{9}), always nonzero, is irreducible precisely when $\Lambda^3U_1=\{0\}$, that is, when $A_n=A_1$, so that $C_{1,2}(U_1)=U_1$.  In this case, the algorithm designates $V_3=U_1\otimes [0,0,2,0,\ldots,0],$ which becomes $U_1\otimes{\mathbb C}$ when $m=2$.  For $m>2$, $$V_1\otimes V_3=(\mbox{ adj }A_1\oplus{\mathbb C})\otimes ([2,0,2,0,\ldots,0]\oplus [1,0,1,1,0,\ldots,0] \oplus [0,0,0,2,0,\ldots,0]),$$  and for $m=2$, 
$$V_1\otimes V_3=(\mbox{ adj }A_1\oplus{\mathbb C})\otimes U_2.$$
This leaves us with 
$$\cT_3^c=\mbox{ adj }A_1\otimes [0,0,0,2,0,\ldots,0]\oplus {\mathbb C}\otimes [1,0,1,1,0,\ldots,0]$$\begin{equation}\label{10}\oplus {\mathbb C}\otimes [0,0,0,2,0,\ldots,0]),\end{equation} when $m>2$, and $\cT_3^c=\{0\}$ when $m=2$.  In the latter case, the algorithm terminates successfully with $j=2$. 

As long as $m>2$, the expression in (\ref{10}) is not irreducible.  It
is equivalent to $\mbox{ adj }A_1\oplus\mbox{ adj }A_3\oplus{\mathbb C}$ if and only if $m=3$.  In this case, the algorithm terminates successfully with $j=4$, $V_4=\mbox{ adj }A_1\oplus\mbox{ adj }A_3$.   If $m>3$, the algorithm aborts with $j=3$.

Next we consider when (\ref{9}) is $\adg\oplus{\mathbb C}$.  We need $C_{1,2}(U_1)=\mbox{ adj }A_2$, $[0,0,2,0,\ldots,0]={\mathbb C}$, and $\Lambda^3U_1={\mathbb C}$.  This happens precisely when $n=m=2$, giving us a successful termination of the algorithm with $j=3$.
 $\Box$ 

\begin{lem}  Let $\mfl_2=A_m$,  and $U_2=[0,1,0,\ldots,0]$. 
The algorithm terminates successfully if and only if one of the following holds. \begin{enumerate}
\item $m=3$, in which case $j=2$ and $\g$ corresponds to $\alpha_{n+1}$-graded $D_{n+4}$; 
\item $1\leq n \leq 4$ and $m=4$; in these cases, $j=2, 3, 5, 5$ respectively and $\g$ corresponds to $\alpha_5$-graded $E_\ell$, $\ell=6,7,8$, or $\alpha_5$-graded $E_8^{(1)}$,  respectively;
\item $n=1$ or $2$, and $m=5$; in these cases, $j=3$ and $\g$ corresponds to $\alpha_3$-graded $E_7$ or $\alpha_3$-graded $E_7^{(1)}$, respectively;
\item $n=1$, and $m=6$ or $7$; in these cases, $j=4$, and $\g$ corresponds to $\alpha_3$-graded $E_8$ or $\alpha_3$-graded $E_8^{(1)}$ respectively.
 \end{enumerate}
In all other cases, the algorithm aborts with $j\leq 4$.
\end{lem}

{\bf Proof}  We have $\Lambda^2 U_2=[1,0,1,0,\ldots,0]$ and $S^2 U_2=[0,2,0,\ldots,0]+[0,0,0,1,0,\ldots,0]$ so that $V_2=[0,1,0,\ldots,0]\otimes [0,0,0,1,0,\ldots,0]$ and $$V_1\otimes V_2=([1,1,0,\ldots,0]+[0,0,1,0,\ldots,0])\otimes([0,1,0,1,0,\ldots,0]+$$$$[1,0,0,0,1,0,\ldots,0]+[0,0,0,0,0,1,0,\ldots,0]).$$ 
This gives us 
$\cT_2^c=[1,1,0,\ldots,0]\otimes[0,0,0,0,0,1,0,\ldots,0]+[0,0,1,0,\ldots,0]\otimes $\begin{equation}\label{n=1}[1,0,0,0,1,0,\ldots,0]+ [0,0,1,0,\ldots,0]\otimes [0,0,0,0,0,1,0,\ldots,0].\end{equation}
When $n=2$ and $m=5$, $\cT_2^c\cong\adg +\mbc$ so the algorithm terminates successfully with $j=2$.  All other combinations of $n$ and $m$ with $n>1$ and $m>4$ yield $V$ inadmissible, the algorithm aborting here with $j=2$.
When $m=3$, $V_2=[0,1,0,\ldots,0]\otimes\mbc$ for arbitrary $n$. In these cases, $V_1\otimes V_2=\cT_2$ and the algorithm terminates successfully with $j=2$. 

Restrict now to $n=1$, $m>3$.  Eq.~(\ref{n=1}) becomes $[1]\otimes [0,0,0,0,0,1,0,\ldots,0]$, which is zero when $m=4$.  Thus, $n=1$ and $m=4$ yield a successful run of the algorithm terminating with $j=2$.  Taking $m=5$, we get $V_3=[1]\otimes \mbc$.
Continuing, $V_1\otimes V_3=([2]+\mbc)\otimes U_2$ so that $\cT_3^c$ is zero and the algorithm terminates successfully with $j=3$.  Letting $m=6$, we get $V_3=[1]\otimes [0,0,0,0,0,1]$ so that $V_1\otimes V_3=([2]+\mbc)\otimes ([0,1,0,0,0,1]+[1,0,0,0,0,0])$ and $V_4=\mbc\otimes[1,0,0,0,0,0]$. Then $V_1\otimes V_4=[1]\otimes([1,1,0,0,0,0]+[0,0,1,0,0,0])$ which implies that $\cT_4^c=\{0\}$, thus, that the algorithm terminates successfully with $j=4$.  If $m=7$, $V_3=[1]\otimes[0,0,0,0,0,1,0]$ and $V_1\otimes V_3=([2]+\mbc)\otimes ([0,1,0,0,0,1,0]+[1,0,0,0,0,0,1]+\mbc)$, giving us $\cT_3^c=[2]\otimes\mbc+\mbc\otimes[1,0,0,0,0,0,1]+\mbc\otimes\mbc$.  Then $V_4\cong\adg$ and the algorithm terminates successfully with $j=4$.  If $n=1$ and $m>7$, the adjoint representation of $\mfl_2$ does not appear in $\cT^c_3$ so the algorithm aborts at $j=3$.

Next take $n=2$.  We need only consider $m=4$. Here $\cT_2^c=\mbc\otimes [1,0,0,0]$  giving us $V_1\otimes V_3=[1,0]\otimes([1,1,0,0,]+[0,0,1,0])$.  Since $\cT_3^c=\{0\}$, the algorithm terminates successfully with $j=3$.

Now reconsider Eq.~(\ref{n=1}) when $n>2$, $m>3$.  Since $V$ is inadmissible when $m\geq 5$, we only need consider $m=4$.  Then $\cT_2^c=[0,0,1,0,\ldots,0]\otimes [1,0,0,0]$ so that  $\cT_3^c=[0,0,0,1,0,\ldots,0]\otimes [0,0,1,0]$.  The algorithm advances to $j=4$ and we have
$$V_1\otimes V_4=([1,0,0,1,0,\ldots,0]+[0,0,0,0,1,0,\ldots])\otimes$$$$([0,1,1,0]+[1,0,0,1]+\mbc).$$ If $n>4$, the algorithm aborts with $j=4$.  When $n=4$, $\cT^c_4=\adg+\mbc$ so the algorithm terminates successfully with $j=5$. 
When $n=3$, the algorithm advances with $V_4=[1,0,0]\otimes\mbc$. Then $V_1\otimes V_4=\cT_4$ so the algorithm terminates successfully with $j=4$. $\Box$

\begin{lem}  Let $\mfl_2=B_m$ or $D_m$ and $U_2=[1,0,\ldots,0]$. The algorithm terminates successfully with $j=2$, $\g$ corresponding respectively to $\alpha_{n+1}$-graded $B_{n+m+1}$ or $D_{n+m+1}$.
\end{lem}

{\bf Proof} If $U_2$ is the standard representation of $B_m$, then $\Lambda^2U_2=\mbox{ adj }B_m$ and $S^2U_2=[2,0,\ldots,0] \oplus{\mathbb C}$. Thus, $$\cT_1^c=\Lambda^2U_1\otimes{\mathbb C},$$ which is always irreducible.  The algorithm designates $V_2=\Lambda^2U_1\otimes{\mathbb C}$ and we have $$V_1\otimes V_2=(C_{1,2}(U_1)\oplus\Lambda^3U_1)\otimes U_2.$$ 
Since $\cT_2^c=\{0\}$, the algorithm terminates successfully with $j=2$.

The argument applies exactly as stated if $U_2$ is the standard representation of $D_m$.
 $\Box$ 

\begin{lem}  Let $\mfl_2=D_5$ and $U_2=[0,0,0,0,1]$. When $n=1, 2$ or $3$, the algorithm terminates successfully, with $j=2$ when $n=1$, and $j=4$ when $n=2$ or $3$. $\g$ corresponds respectively to $\alpha_6$-graded $E_7$, $\alpha_6$-graded $E_8$, or $\alpha_6$-graded $E_8^{(1)}$. Otherwise, the algorithm aborts with $j=3$.  \end{lem}

{\bf Proof}  Under our hypotheses, $S^2U_2=[0,0,0,0,2]\oplus[1,0,0,0,0]$
so that $\cT_1^c$ is always irreducible and the algorithm always advances to $j=2$.  This gets us $V_2=\cT_1^c=\Lambda^2U_1\otimes[1,0,0,0,0]$ and  $$V_1\otimes V_2=(C_{1,2}(U_1)\oplus\Lambda^3U_1)\otimes([1,0,0,0,1] \oplus [0,0,0,1,0]).$$  Note that $$\cT_2^c=\Lambda^3U_1\otimes [0,0,0,1,0]$$  is zero if $n=1$, irreducible otherwise.  We see, then, that when $n=1$, the algorithm terminates successfully with $j=2$. If $n>1$, the algorithm advances to $j=3$ with $V_3=\Lambda^3U_1\otimes [0,0,0,1,0].$  This gives us 
$$V_1\otimes V_3=(C_{1,3}(U_1)\oplus \Lambda^4U_1)\otimes ([0,0,0,1,1] \oplus \mbox{ adj }D_5 \oplus{\mathbb C}),$$  so that \begin{equation}\label{11}\cT_3^c=C_{1,3}(U_1)\otimes{\mathbb C}\oplus \Lambda^4U_1\otimes\mbox{ adj }D_5\oplus \Lambda^4U_1\otimes {\mathbb C}.\end{equation}
The expression in (\ref{11}) is irreducible if and only if $\Lambda^4U_1=\{0\}$, that is, if and only if $A_n=A_2$, in which case $V_4=U_1$.  Continuing, we have $$V_1\otimes V_4=(S^2 V_1\oplus \Lambda^2 V_1)\otimes V_2$$ so that $\cT_4^c=\{0\}$, and the algorithm terminates successfully with $j=4$.

Next consider that the expression in (\ref{11}) is equivalent to $\adg\oplus {\mathbb C}$ if and only if $\Lambda^4V_1={\mathbb C}$, which is true if and only if $A_n=A_3$.  In this case, the algorithm terminates successfully with $j=4$.

For all values of $n>3$, the expression in (\ref{11}) forces the algorithm to abort with $j=3$.$\Box$ 

\begin{lem} Let $\mfl_2=E_6$, $U_2=[1,0,0,0,0,0]$.  The algorithm terminates successfully if and only if $n=1$ or $2$, in which case $j=3$ and $\g$ corresponds respectively to $\alpha_7$-graded $E_8$ or $\alpha_7$-graded $E_8^{(1)}$. Otherwise, the algorithm aborts with $j=2$. \end{lem}

{\bf Proof}  Under the hypotheses, say $U_2=[1,0,0,0,0,0]$. We have
$$\cT_1^c=\Lambda^2U_1\otimes [0,0,0,0,0,1]$$ so that the algorithm always advances to $j=2$.  Then
$$U_1\otimes U_2=(C_{1,2}(U_1)\oplus\Lambda^3U_1)\otimes([1,0,0,0,0,1] \oplus \mbox{ adj }E_6 \oplus{\mathbb C}),$$ leaving us with \begin{equation}\label{12}\cT_2^c=C_{1,2}(U_1)\otimes{\mathbb C}\oplus\Lambda^3U_1\otimes\mbox{ adj }E_6\oplus\Lambda^3U_1\otimes{\mathbb C}.\end{equation} This is irreducible if and only if $\Lambda^3U_1=\{0\}$, that is, if and only if $A_n=A_1$.  In this case, $V_4=U_1\otimes{\mathbb C}$ and we continue, to find $V_1\otimes V_4=\cT_4.$  This is a successful termination of the algorithm with $j=3$. 

The expression in (\ref{12}) is equivalent to $\adg\oplus{\mathbb C}$ if and only if $\Lambda^3U_1={\mathbb C}$, that is, if and only if $A_n=A_2$.  In this case, the algorithm terminates successfully with $j=3$. 

In all other cases, the expression in (\ref{12}) forces the algorithm to abort with $j=2$. $\Box$ 

This exhausts the cases where $\cT_1^c$ irreducible.  
Next, we consider conditions on $V=U_1\otimes U_2$ that force $\cT_1^c\cong\adg+\mbc.$  Since $$\Lambda^2 V=\Lambda^2 U_1\otimes S^2 U_2\oplus S^2U_1\otimes\Lambda^2U_2,$$ each summand of $\Lambda^2 V$ would have the form $W+\mbox{ ad }\mfl_i\otimes \mbc$, where $W$ is itself a tensor product, either of something with $\mbox{ ad }\mfl_i$ or something with $\mbc$.  Note in particular that neither $\Lambda^2 U_i$ nor $S^2 U_i$ can have the form $\mbox{ ad }\mfl_i\oplus \mbc$.  Note further that $S^2 U_i\neq \mbc$ and that $\Lambda^2 U_i=\mbc$ if and only if $U_i$ is the standard representation of $A_1$, in which case $S^2 U_i=\mbox{ ad }A_1$.   These allow us to limit the criteria determining $U_1$ and $U_2$ to the following.
\begin{enumerate}
\item $\Lambda^2 U_i=\mbox{ad }\mfl_i$, $S^2 U_i=Z+\mbc$, where $Z$ is irreducible; 
\item $U_1$ is the standard representation of $A_1$, $S^2 U_2=Z+\mbox{ ad }\mfl_2$, and $\Lambda^2 U_2=W+\mbc$, where $Z$ and $W$ are irreducible.
\end{enumerate}
What follows is an complete list of generalized cominuscule representations $U$ that satisfy (1). (LiE \cite{lie} helps verify that the list is exhaustive.) Taking $U_1$, $U_2$ to be any pair from the list, we get $\cT_1^c\cong\adg+\mbc$.  

\begin{enumerate}
\item $\mfl=A_1$, $U=[2]$;
\item $\mfl=A_3$, $U=[0,1,0]$; 
\item $\mfl=B_n$, $U=[1,0,\ldots,0]$;
\item $\mfl=D_n$, $U=[1,0,\ldots,0]$.
\end{enumerate}

There are sixteen pairs we can choose.  With $V=U_1\otimes U_2$, the algorithm terminates successfully with $j=2$ and $\g$ corresponding to one of the following types of lsn-graded affine algebras:  $\alpha_k$-graded $D_{n+1}^{(2)}$, $k\in\{1,\ldots, n-1\}$; $\alpha_k$-graded $D_n^{(1)}$, $k\in\{3,\ldots, n-3\}$; 
$\alpha_k$-graded $B_n^{(1)}$, $k\in\{3,\ldots,n-1\}$.  Further, any such lsn-graded algebra can be produced by the algorithm for some choice of $U_1$ and $U_2$ on this list.

Next, consider the following exhaustive list of generalized cominuscule representations $U$ with $S^2 U=Z+\mbox{ ad }\mfl$, and $\Lambda^2 U=W+\mbc$, where $Z$ and $W$ are irreducible.  
\begin{enumerate}
\item $\mfl=A_1$, $U=[3]$;
\item $\mfl=A_5$, $U=[0,0,1,0,0]$;
\item $\mfl=C_3$, $U=[0,0,1]$;
\item $\mfl=D_6$, $U=[0,0,0,0,1,0]$;
\item $\mfl=E_7$, $U=[0,0,0,0,0,0,1]$.
\end{enumerate}

Taking $V=U_1\otimes U_2$, with $U_1$ the standard representation of $A_1$, and $U_2$ from this list, we get $\g$ corresponding respectively to:  $\alpha_2$-graded $G_2^{(1)}$; $\alpha_3$-graded $E_6^{(1)}$; $\alpha_1$-graded $F_4^{(1)}$; $\alpha_1$-graded $E_7^{(1)}$; $\alpha_8$-graded $E_8^{(1)}$.

The following summarizes our accounting in case $V=U_1\otimes U_2$. We omit cases that terminate successfully with $j=1$.

\begin{prop}
Let $U_1$ be a generalized cominuscule representation for finite type $\mfl_1$, $U_2$ a generalized cominuscule representation for finite type $\mfl_2$.  The algorithm either aborts with $j\leq 4$ or it terminates successfully with $1<j\leq 4$ and $\g$ corresponding to one of the following types of lsn-graded algebras:  $\alpha_k$-graded $B_n$, $1<k<n-1$; $\alpha_k$-graded $B_n^{(1)}$, $2<k<n-1$; $\alpha_k$-graded $D_n$, $1<k<n-2$; $\alpha_k$-graded $D_n^{(1)}$, $2<k<n-2$; $\alpha_k$-graded $D_{n+1}^{(2)}$, $0<k<n$; $\alpha_k$-graded $E_n$, $\alpha_k$ graded $E_n^{(1)}$, $k=3$ or $5\leq k\leq n-1$; $\alpha_2$-graded $E_6^{(1)}$; $\alpha_1$-graded $E_7^{(1)}$; $\alpha_8$-graded $E_8^{(1)}$; $\alpha_3$-graded $E_6^{(2)}$; $\alpha_2$-graded $F_4$; $\alpha_k$-graded $F_4^{(1)}$, $k=1,2$; $\alpha_2$-graded $G_2^{(1)}$. $\Box$\end{prop}

\subsection{$\gk$ Simple}
Here we detail the cases where $\cT_1^c$ is irreducible. Since $\gk$ is simple, the irreducible cases subsume the cases where $\cT_1^c\cong\adg$.

$\cT_1^c$ is zero precisely when $\Lambda^2 V$ is irreducible. 
We list of these cases to have a complete account.

\begin{enumerate}
\item $\gk=A_n$, $V=[1,0,\ldots,0]$, $\g$ corresponds to $\alpha_1$-graded $A_{n+1}$;
\item $\gk=A_n$, $V=[2,0,\ldots, 0]$, $\g$ corresponds to $\alpha_{n}$-graded $C_{n+1}$;
\item $\gk=A_n$, $V=[0,1,0,\ldots,0]$, $\g$ corresponds to $\alpha_n$-graded $D_{n+1}$;
\item $\gk=B_n$ or $D_n$, $V=[1,0,\ldots,0]$, $\g$ corresponds to a type $\alpha_1$-graded $B_{n+1}$ or $\alpha_1$-graded $D_{n+1}$ algebra respectively.
\end{enumerate}

$\cT^c_1$ is irreducible precisely when $\Lambda^2V$ has exactly two irreducible components.  We list the representations $V$ with this property.   In cases where $\cT^c_1=\adg$, the algorithm terminates successfully with $j=2$.  In cases where $\cT^c_1={\mathbb C}$, $V_1\otimes V_2$ is irreducible, also resulting in a successful termination of the algorithm at $j=2$.  We note these in the list and analyze the rest of the cases through the lemmas that follow.

\begin{enumerate}
\item $\gk=A_n$, $n\geq 5$, $V_1=[0,0,1,0,\ldots,0]$; $V_2=\cT_1^c=[0,0,0,0,0,1,0,\ldots,0]$;  when $n=5$, $V_2=\mbc$ and the algorithm terminates with $j=2$, $\g$ corresponding to $\alpha_2$-graded $E_6$.  
\item $\gk=A_n$, $n\geq 7$, $V_1=[0,0,0,1,0,\dots,0]$, $\cT_1^c=[1,0,0,0,0,0,1,0,\ldots,0]$; when $n=7$, the algorithm terminates successfully with $j=2$ and $\g$ corresponding to $\alpha_2$-graded $E_7^{(1)}$.
\item $\gk=A_n$, $V_1=[3,0,\ldots,0]$ or $[4,0,\ldots,0]$;  $\cT_1^c=[0,3,0,\ldots,0]$ or $[2,3,0\ldots,0]$ respectively.
\item $\gk=A_n$, $V_1=[0,2,0,\ldots,0]$, $\cT_1^c=[1,0,1,1,0,\ldots,0]$; when $n=3$, the algorithm terminates with $j=2$, $\g$ corresponding to $\alpha_3$-graded $A_{5}^{(2)}$.
\item $\gk=B_n$ or $D_n$, $V_1=[2,0,\ldots,0]$, $\cT_1^c=\mbox{ adj }B_n$, respectively, 
$\cT_1^c=\mbox{ adj }D_n$; the algorithm terminates successfully with $\g$ corresponding to $\alpha_n$-graded $A_{2n}^{(2)}$, respectively $\alpha_n$-graded $A_{2n-1}^{(2)}$.
\item $\gk=C_3$ or $C_4$, $V_1=[0,\ldots,0,1]$, $\cT_1^c=\mathbb C$, respectively $\cT^c_1=\adg$;  the algorithm terminates successfully with $j=2$, $\g$ corresponding respectively to $\alpha_2$-graded $F_4$, or $\alpha_4$-graded $E_6^{(2)}$.
\item $\gk=E_6$, $V_1=[2,0,0,0,0,0]$, $\cT_1^c=[0,0,1,0,0,1]$.
\item $\gk=E_7$, $V_1$ the cominuscule representation, $\cT_1^c={\mathbb C}$; the algorithm terminates successfully with $j=2$ and $\g$ corresponding to $\alpha_8$-graded $E_8$. \end{enumerate}

\begin{lem} Let $\gk=A_n$, $n>5$, $V_1=[0,0,1,0,\ldots,0]$. Unless $n=6,7$ or $8$, $V$ is inadmissible and the algorithm aborts with $j=2$.  When $n=6$ or $7$, the algorithm terminates successfully with $j=2$, and $\g$ corresponding respectively to $\alpha_2$-graded $E_7$ or $\alpha_2$-graded $E_8$.  When $n=8$, the algorithm terminates successfully with $j=2$ and $\g$ corresponding to $E_8^{(1)}$.
\end{lem}

{\bf Proof} When $n=6$, $\cT_1^c=[0,0,0,0,0,1]$, $V_1\otimes V_2=\cT_2$, so the algorithm terminates successfully with $j=2$.  For $n\geq 7$, we have $$V_1\otimes V_2=[0,0,1,0,0,1,0,\ldots,0]\oplus [0,1,0,0,0,0,1,0,\ldots,0]\oplus $$$$[1,0,0,0,0,0,0,1,0,\ldots,0]\oplus [0,0,0,0,0,0,0,0,1,0,\ldots,0].$$ The algorithm thus aborts when $n>8$.  If $n=7$, $V_3= [1,0,0,0,0,0,0]$.  Then $$V_1\otimes V_3= [1,0,1,0,0,0,0]\oplus [0,0,0,1,0,0,0],$$ so that $\cT_3^c=\{0\}$.  The algorithm thus terminates successfully with $j=2$. When $n=8$, $V_3=\adg$ so the algorithm terminates successfully with $j=3$. $\Box$

\begin{lem} If $\gk=A_n$, $n>7$, and $V_1=[0,0,0,1,0,\dots,0]$, the algorithm aborts with $j=2$\end{lem}

{\bf Proof}  We have $V_2=[1,0,0,0,0,0,1,0,\ldots,0]$ so that when $n>7$, $V_1\otimes V_2$ has at least four nontrivial components.  $\Box$

\begin{lem} Let $\gk=A_n$ and $V_1=[3,0,\ldots,0]$.  The algorithm terminates successfully if and only if  $n=1$ or $2$.  When $n=1$, $j=2$ and $\g$ corresponds to $\alpha_2$-graded $G_2$.  When $n=2$, $j=3$ and $\g$ corresponds to $\alpha_2$-graded $D_4^{(3)}$. 

If $\gk=A_n$ and $V_1=[4,0,\ldots,0]$, the algorithm aborts unless $n=1$.  In this case, the algorithm terminates successfully with $j=2$ and $\g$ corresponding to $\alpha_1$-graded $A_2^{(2)}$.  
\end{lem}

{\bf Proof}  If $V_1=[3,0,\ldots,0]$, then $V_2=[0,3,0,\ldots,0]$ so that $$V_1\otimes V_2=[3,3,0,\ldots,0] \oplus [2,2,1,0,\ldots,0] \oplus [1,1,2,0,\ldots,0] \oplus [0,0,3,0,\ldots,0].$$  It is clear that the algorithm aborts with $j=2$ when $n>2$.  When $n=2$, $V_3=\adg$, and the algorithm terminates successfully with $j=3$.  When $n=1$, $V_2=\mbc$ so $V_1\otimes V_2=\cT_2$, giving us a successful termination of the algorithm with $j=2$.  

When $k=4$, $V_2=[2,3,0,\ldots,0]$ so that when $n=1$, $V_2=\adg$ and the algorithm terminates successfully with $j=2$.  Otherwise, $V_1\otimes V_2$ has eleven nontrivial components and the algorithm aborts. $\Box$ 

\begin{lem} If $\gk=A_n$, $n>3$, and $V=[0,2,0,\ldots,0]$, then $V$ is inadmissible, and the algorithm aborts with $j=2$. \end{lem}

{\bf Proof} We have $V_2=\cT_1^c=[1,0,1,1,0,\ldots,0]$ so that
$V_1\otimes V_2$ has more than five nontrivial irreducible components. 
$\Box$ 

\begin{lem} If $\gk=E_6$ with $V_1=[2,0,0,0,0,0]$, the algorithm aborts with $j=2$. \end{lem}

{\bf Proof} $V_2=[0,0,1,0,0,1]$ so that $V_1\otimes V_2$ has more than four nontrivial components, forcing the algorithm to abort.  $\Box$ 

 This completes the proof that there is an admissible representation associated to any affine or finite Lie algebra with lsn-grading and further,  that these are the only admissible representations.   In the next section, we verify that the algorithms actually produce Lie algebras with finite or affine type root systems. 

\section{Proof of the Algorithm, Part II: The Structure of $\g$}
Here we establish how the structure of an affine or finite type Lie algebra can be extended from $\gk$ to a $\mbz$-graded vector space $\g$ produced from the Affine or Finite Algorithm.  

Suppose $\g$ is constructed via the Affine Algorithm.  Let $\gamma_1=\sum_{i\neq k} c_i\alpha_i$ be the lowest weight of $\g_1$. Define $a_k$ to be the least common denominator of the $c_i$s.  Since weights of $\g_1$ differ by elements of $Q_+$, $a_k$ is the least common denominator when we write any weight of $\g_1$ as a linear combination of elements of $\Pik$. Take $r=j/a_k$, where $j$ is minimal with $\g_j\cong \gk$.  We claim that $j$ is a multiple of $a_k$:  since $\g_j\cong \gk$, each $\alpha_i\in\Pik$ must be a weight of $\g_j$.  The highest weight of $\g_j$ has the form $j\lambda_1$ less $2j-1$ positive roots of $\gk$.  In particular, $j\lambda_1$ is in the root lattice of $\gk$ so $j/a_k$ must be a positive integer as claimed.

If the Finite Algorithm produces $\g$, let $t_k=1/(2-(\lambda_1,\lambda_1))$ where $\lambda_1$ is the highest weight of $\g_1$. 

Let $B$ be the Killing form on $\gk$.   Use $(.,.)$ to denote its rescaling as per \S4.4, that is, the long roots associated to a simple component of $\gk$ have length $2/n$ where $n$ is the degree of the representation $\g_1$ on that component of $\gk$.  

The Affine Algorithm defines $\g_0=\gk\oplus \mbc+\mbc$.  Identify one copy of $\mbc$ with $\mbc K$ and the other copy with $\mbc d_k$, where $K.x=0$ for all $x\in \g$ and $d_k.x_t=tx$ if $x\in\g_t$. In particular, define $[K\;x]=[d_k\;x]=0$ for all $x\in\g_0$ so $\g_0$ is a Lie algebra and $\g$ is a $\g_0$ module. Use Lemma~\ref{affinefactslem} to extend the definition of  $(.,.)|_{\hk}$ to all of $\ho:=\hk+\mbc K+\mbc d_k$.  Decree that $K$ and $d_k$ are orthogonal to $x\in\g$ if $x\not\in\ho$. Let $\nu:\ho\rightarrow \hod$ be the isomorphism determined by the nondegenerate form on $\ho$. Define $\delta:=\nu(K)$,  $\Lambda_k:=\nu (d_k)$, and $\alpha_k:= \gamma_1+(1/a_k)\delta=-\lambda_1+(1/a_k)\delta$.  This puts $\alpha_k$ in $\hod$ with $(\alpha_k,\alpha_i)=(\gamma_1,\alpha_i)=-(\lambda_{-1},\alpha_i)$ for all $\alpha_i\in\Pik$.  Note that $\|\alpha_k\|^2=2$. It is easy to verify that the net effect of all this on the roots of $\gk$ is to append a node representing $\alpha_k$ to $\Dk$ as it represents the module $\g_{-1}$:  $\alpha_k$ is long and shares $n$ edges with any node of $\Dk$ labeled $n$.   Identify $\g_{\alpha_k}$ with the $\gamma_1$ weight space of $\g_1$ and  $\g_{-\alpha_k}$ with the $\lambda_{-1}$ weight space of $\g_{-1}$.  

The Finite Algorithm defines $\g_0=\gk\oplus \mbc$.  Identify $\mbc$ with $\mbc d_k$ where $d_k.x=tx$ for $x\in\g_t$.  Use Lemma~\ref{finitefactslem} to extend the definition of $(.,.)|_{\hk}$ to all of $\ho:=\hk+\mbc d_k$.  Decree that $d_k$ is orthogonal to $x\in\g$ if $x\not\in\ho$.  Let $\nu:\ho\rightarrow \hod$ be the isomorphism determined by the nondegenerate form on $\ho$. Define $\Lambda_k:=\nu (d_k)$, and $\alpha_k:= \gamma_1+(1/t_k)\Lambda_k=-\lambda_{-1}+(1/t_k)\Lambda_k$. The rest of the previous paragraph now applies without change.

We now have the elements necessary to define all the brackets on the graded pieces of $\g$ as in \S3.  The remainder of the discussion goes towards verifying that the $\alpha_k$-extended $\Dk$ actually describes $\g$ as a Lie algebra with those brackets.

\subsection{Brackets}
   Let $[\cdot,\cdot]:\g_{-1}\otimes\g_1\rightarrow\g_0$ be given by 
\begin{equation}\label{13} [u_{-1}\;u]=-\sum_i(u_{-1},X_i.u)Y_i+\kappa(u_{-1},u)\end{equation} where in the affine case, $\kappa=(-1/a_k)K$ and in the finite case, $\kappa=(-1/t_k)d_k$. 
Defining the rest of the brackets as in \S3, we must
 show that 
$\cT_1$ is in the kernel of the bracket defined on $\Lambda^2 V=\Lambda^2\g_1$ and, for $i>2$, that $\cT_i$ is in the kernel of the bracket defined on $\g_1\otimes\g_i$.  We verify that the bracket maps $\g_1\otimes\g_i$ onto $\g_{i+1}$ by producing a highest weight vector of $\g_{i+1}$ in terms of the bracket.  Operators help  advance these arguments.

Let $\{X_i\}$ and $\{Y_i\}$ be Killing dual bases of $\gk$. Define $\Phi$ on $\Lambda^2\g_1$ by $\Phi(u\wedge v)=\sum_i X_i.v\wedge Y_i.u$.  

\begin{prop} \label{phiprop}
In the affine case, $\Phi|_{\cT_1}\equiv 0.$  In the finite case, $\Phi|_{\cT_1}\equiv -1/t_k.$ 
\end{prop}

{\bf Proof}  We have $C.(u\wedge v)=2c_1(u\wedge v)+2\Phi(u\wedge v).$  
The result follows by Proposition~\ref{casprop}. $\Box$ 

We take $[u\;v]=\Phi(v\wedge u)$ in the affine case and in the finite case, $[u\;v]=\Phi(v\wedge u)+\F{1}{t_k}v\wedge u$  (cf. Eqs.~(3) and (4), \S3.3.)  Invoking Proposition~\ref{phiprop} we get the following.

\begin{thm}  The bracket is identically zero on $\cT_1$. $\Box$\end{thm}

\begin{lem}\label{hwvlem}  A highest weight vector of $\g_2$ has the form $[v^+\;f_{\beta_1}f_{\beta_2}f_{\beta_3}.v^+]$, the $\beta_i$s as given in Proposition~\ref{subwtprop}. \end{lem}

{\bf Proof}  In the affine case, we have
$$[v^+\;f_{\beta_1}f_{\beta_2}f_{\beta_3}.v^+]=-\sum_i X_i.v^+\wedge Y_i.f_{\beta_1}f_{\beta_2}f_{\beta_3}.v^+$$$$= -(\lambda_1,\lambda_1-\beta_1-\beta_2-\beta_3)v^+\wedge f_{\beta_1}f_{\beta_2}f_{\beta_3}.v^++\sum v_\nu\wedge v_\mu$$$$=v^+\wedge f_{\beta_1}f_{\beta_2}f_{\beta_3}.v^+ +\sum v_\nu\wedge v_\mu$$
where $v_\nu$ and $v_\mu$ are weight vectors in $\g_1$ with $\nu<\lambda_1$ and $\mu>\lambda_1-\beta_1-\beta_2-\beta_3$.  Note in particular that $\nu$ is always of the form $\lambda_1-\alpha$ for a positive root $\alpha$ and $\mu=\lambda_1-\beta_1-\beta_2-\beta_3+\alpha$. Since $v^+\wedge f_{\beta_1}f_{\beta_2}f_{\beta_3}.v^+\neq 0$,  $[v^+\;f_{\beta_1}f_{\beta_2}f_{\beta_3}.v^+]\neq 0$.  As the component of $[v^+\; f_{\beta_1}f_{\beta_2}f_{\beta_3}.v^+]$ in $\cT_1$ is zero, $[v^+\; f_{\beta_1}f_{\beta_2}f_{\beta_3}.v^+]$ belongs to $\g_2$ and by weight considerations, the result follows.

In the finite case, 
$$[v^+\;f_{\beta_1}f_{\beta_2}f_{\beta_3}.v^+]=-\sum_i X_i.v^+\wedge Y_i.f_{\beta_1}f_{\beta_2}f_{\beta_3}.v^+-\F{1}{t_k}v^+\wedge f_{\beta_1}f_{\beta_2}f_{\beta_3}.v^+=$$$$-\left((\lambda_1,\lambda_1-\beta_1-\beta_2-\beta_3)+\F{1}{t_k}\right)v^+\wedge f_{\beta_1}f_{\beta_2}f_{\beta_3}.v^+ +\sum v_\nu\wedge v_\mu=$$$$\left(1-\F{1}{t_k}\right)v^+\wedge f_{\beta_1}f_{\beta_2}f_{\beta_3}.v^+ +\sum v_\nu\wedge v_\mu.$$
For reasons cited above, this is nonzero and in $\g_2$ when $t_k\neq 1$.  Even if $t_k=1$, though, it is clear that $f_{\beta_1}.v^+\wedge f_{\beta_2}f_{\beta_3}.v^+$, a multiple of one of the $v_\nu\wedge v_\mu$ terms, is nonzero.  The constant factor itself would be a nonzero multiple of $(\lambda_1-\beta_2-\beta_3,\beta_1)=1$.  Thus in these cases as well, the result holds.  $\Box$ 

\begin{lem}
For $u_{-1}\in\g_{-1}$ and $v,\;w\in\g_1$, we have, in the affine case
$$[u_{-1}\;[v\;w]\;]=\sum_i - (u_{-1},X_i.v)Y_i.w+(u_{-1},X_i.w)Y_i.v.$$
In the finite case,
$$[u_{-1}\;[v\;w]\;]=\sum_i - (u_{-1},X_i.v)Y_i.w+(u_{-1},X_i.w)Y_i.v
-\F{1}{t_k}(u_{-1},v)w+\F{1}{t_k}(u_{-1},w)v.$$

\end{lem}

{\bf Proof} By the Jacobi identity we have
$$[u_{-1}[v\;w]\;]=[\;[u_{-1}\;v]w]+[v[\;u_{-1}\;w]\;]=[\;[u_{-1}\;v]w]-[\;[u_{-1}\;w]v].$$ The result then follows application of Eqs.~(3), (4) from \S3.3, along with Eq.~(\ref{13}).
$\Box$

For $t<j$, define $\Psi$ on $\g_1\otimes\g_t$ by 
$$\Psi(u\otimes[u_1\;\ldots\;u_t])=\sum_i X_i.u\otimes Y_i.[u_1\;\ldots\;u_t]+X_i.u_2\otimes [u_1\;Y_i.[u_3\;\ldots\;u_t]\;]+$$$$X_i.u_3\otimes [u_1\;u_2\;Y_i.[u_4\;\ldots\;u_t\;]\;]+\ldots +X_i.u_{t-1}\otimes[u_1\;\ldots\;u_{t-2}\;Y_i.u_t\;]$$$$-X_i.u_t\otimes [u_1\;\ldots\;u_{t-2}\;Y_i.u_{t-1}\;].$$

\begin{prop}  $\Psi$ is $\gk$-equivariant.
\end{prop}

{\bf Proof}  We can view $\Psi$ as a sum of compositions of two kinds of mappings:  (1) permutations on multivectors in $\g_1\otimes\g_t\subset \g_1\otimes\ldots\otimes\g_1$; and (2) Casimir polarizations on multivectors, that is, mappings of the form 
$$u\otimes u_1\otimes\ldots\otimes u_t\mapsto \sum_i u\otimes\ldots \otimes u_{\ell-1}\otimes X_i.u_\ell\otimes u_{\ell+1}\otimes \ldots\otimes u_{m-1}\otimes Y_i.u_m\otimes u_{m+1}\otimes \ldots u_t.$$
This is a tedious but straightforward verification that one enacts starting with the observation that
$$\Psi(u_1\otimes u_2)=\Phi(u_1\wedge u_2)=\sum_i X_i.u_1\otimes Y_i.u_2-X_i.u_2\otimes Y_i.u_1$$
can be realized as a composition of mappings
$$u_1\otimes u_2\mapsto u_1\otimes u_2- u_2\otimes u_1$$
and 
$$u_1\otimes u_2\mapsto \sum_i X_i.u_1\otimes Y_i.u_2.$$
Both types of mappings are $\gk$-equivariant.  The first is the difference between the identity and the usual action of the symmetric group on $\g_1\otimes \g_1$.  The second is $\F{C}{2}-c_1$. Since $C$ is $\gk$-equivariant, its polarization is as well.  
 $\Box$ 

It is convenient to think of $\Psi$ as a sum of three operators so define
$$\Psi_1(u\otimes [u_1\;\ldots\;u_t\;])=\sum_i X_i.u\otimes Y_i.[u_1\;\ldots\;u_t\;],$$
$$\Psi_2(u\otimes [u_1\;\ldots\;u_t\;])=\sum_i X_i.u_1\otimes [u\;Y_i.[\;u_2\;\ldots\;u_t\;]\;],$$
$$\Psi_3=\Psi-\Psi_1-\Psi_2.$$

The argument we advanced to support the claim that $\Psi$ is $\gk$-equivariant applies to give us $\Psi_1$, $\Psi_2$ and $\Psi_3$ $\gk$-equivariant as well.

Recall that $Y_{1,t}$ is the Cartan product of $\g_1$ and $\g_t$.  Let $U_{1,t}$ be the submodule of $\g_1\otimes\g_t$ with highest weights of the form $\lambda_1+\lambda_t-\alpha$, $\alpha\in\Deltak_+$.

We note the following.

\begin{lem}\label{psi1lem} In the affine case, $\Psi_1|_{Y_{1,t}}\equiv 1$ and ${\Psi_1}|_{U_{1,t}}\equiv 1-(\rho,\alpha)$.  In the finite case, ${\Psi_1}|_{Y_{1,t}}\equiv 1-\F{t}{t_k}$ and ${\Psi_1}|_{U_{1,t}}\equiv 1-\F{t}{t_k}-(\rho,\alpha)$.\end{lem}

{\bf Proof} It is clear that $\ds \Psi_1=\F{C-(c_1+c_t)}{2}$. The result then follows Corollary~\ref{lam1lamt}. $\Box$ 

Let $w^+=\sum v_\nu\otimes v_\mu$ be a highest weight vector of $\g_t$, so that $v_\nu\in\g_1$ and $v_\mu\in\g_{t-1}$.  Take $v_\nu\otimes v_\mu$ so that $\nu$ is highest among the weights of $\g_1$ associated to $w^+$.  

Since $e_\alpha.w^+=0$ for all positive root vectors $e_\alpha\in\gk$, it must be the case that $(e_\alpha.v_\nu)\otimes v_\mu=0$ for all $\alpha$ or that $w^+$ has a component of the form $v_{\nu+\alpha}\otimes v_{\mu-\alpha}$.  The latter contradicts our choice of $\nu$ as maximal so $e_\alpha.v_\nu\otimes v_\mu=0$ for all $\alpha\in\Deltak_+$ implies $w^+$ has a component of the form $v^+\otimes v_\mu$. 

\begin{prop}\label{mystery}
If $s\leq j$, a highest weight vector for $\g_s$ has the form
\begin{equation}\label{mysteq} {v_s}^+=\overbrace{[v^+\;f_{\gamma_1}f_{\gamma_2}.v^+\;\ldots\;f_{\nu_1}f_{\nu_2}.v^+\;f_{\beta_1}f_{\beta_2}f_{\beta_2}.v^+\;]}^{\mbox{$s$ terms }}\end{equation}
where $\{\gamma_1,\gamma_2\},\ldots,\{\nu_1,\nu_2\},\{\beta_1,\beta_2,\beta_3\}$ are sets of positive roots as in Proposition~\ref{subwtprop}.  Moreover, for $v_{\nu}\in\g_1$, and $v_{\mu}\in\g_{t-1}$, $[v_\nu\;v_\mu]=0$ when $\nu+\mu>\lambda_s$.  In particular, the bracket projects $v_\nu\otimes v_\mu$ into $\g_s$.
\end{prop}

We prove the result by induction.  The case $s=2$ was done in Lemma~\ref{hwvlem}. 
The next sequence of lemmas follow the assumption that the proposition is true for all $s\leq t$.  At the end of those lemmas, we will have a proof of the proposition.

\begin{lem} \label{psi2lem} In the affine case, $\ds\Psi_2|_{Y_{1,t}}\equiv -1$. In the finite case, $\ds\Psi_2|_{Y_{1,t}}\equiv -1-\F{t-1}{t_k}$.
\end{lem}

{\bf Proof}   A highest weight vector of $Y_{1,t}$ has the form 
$$v^+\otimes [v^+\;f_{\gamma_1}f_{\gamma_2}.v^+\;\ldots\;f_{\nu_1}f_{\nu_2}.v^+\;f_{\beta_1}f_{\beta_2}f_{\beta_2}.v^+\;]$$
 which we write $v^+\otimes [\;v^+\;u\;]$.  Note that $u\in\g_{t-1}$.

We have 
\begin{equation}\label{14}\Psi_2(v^+\otimes[v^+\;u])=\sum_iX_i.v^+\otimes [v^+\;Y_i.u].\end{equation}
Terms of (\ref{14}) associated to $X_i$ from positive root spaces are zero because $X_i.v^+=0$.  Terms associated to $Y_i$ from positive root spaces are also zero because then $[v^+\;Y_i.u]$ has a higher weight than the highest weight of $\g_t$.  Thus the nonzero terms in (\ref{14}) are contributed from the part of $\{X_i\}$ and $\{Y_i\}$ that comprise dual bases of $\hk$. Invoking Corollary~\ref{lam1lamt}, we get, in the affine case,
$$\Psi_2(v^+\otimes [v^+\;u])=(\lambda_1,\lambda_{t}-\lambda_1)v^+\otimes [v^+\;u]=-v^+\otimes [v^+\;u].$$
In the finite case, we have
$$\Psi_2(v^+\otimes [v^+\;u])=(\lambda_1,\lambda_{t}-\lambda_1)v^+\otimes [v^+\;u]=\left(-1-\F{t-1}{t_k}\right)v^+\otimes [v^+\;u]. \Box $$

\begin{lem} \label{affpsi} In the affine case, $\Psi|_{Y_{1,t}}\equiv 0.$ \end{lem}

{\bf Proof}  We have $(\Psi_1+\Psi_2)|_{Y_{1,t}}\equiv 0$ in the affine case so it remains to show that $\Psi_3|_{Y_{1,t}}\equiv 0$.  

Let $v^+_{t}=[v^+\;f_{\gamma_1}f_{\gamma_2}.v^+\;\ldots\;f_{\nu_1}f_{\nu_2}.v^+\;f_{\beta_1}f_{\beta_2}f_{\beta_2}.v^+\;]$ be a highest weight vector for $\g_t$. We have 
$$\Psi_3(v^+\otimes v_t^+)=\sum_i X_i.f_{\gamma_1}f_{\gamma_2}.v^+\otimes [v^+\;v^+\;Y_i.[f_{\xi_1}f_{\xi_2}.v^+\ldots \;f_{\nu_1}f_{\nu_2}.v^+\;f_{\beta_1}f_{\beta_2}f_{\beta_2}.v^+\;]]+\ldots$$\begin{equation}\label{15}-X_i.f_{\beta_1}f_{\beta_2}f_{\beta_2}.v^+\otimes[v^+\;v^+\;f_{\gamma_1}f_{\gamma_2}.v^+\ldots\;Y_i.f_{\nu_1}f_{\nu_2}.v^+].\end{equation}
According to Proposition~\ref{subwtprop}, $\lambda_{t-1}$ has the form $(t-1)\lambda_1$ less $2t-3$ positive roots.
Consider that terms of $\g_{t-1}$ appearing in (\ref{15}) with the form
$$[v^+\;v^+\;f_{\gamma_1}f_{\gamma_2}.v^+\ldots Y_i.[f_{\upsilon_1}f_{\upsilon_2}.v^+\ldots \;f_{\nu_1}f_{\nu_2}.v^+\;f_{\beta_1}f_{\beta_2}f_{\beta_2}.v^+\;]\;]$$
are brackets of $t-1$ vectors from $\g_1$.  The lowest weight that could be associated to such a term has the form $(t-1)\lambda_1$ less $2(t-3)+1+1=2t-4$ positive roots.  Any such weight is strictly greater that $\lambda_{t-1}$.  By assumption, these brackets must all be zero, which proves the lemma.
$\Box$ 

Comparing to Eq.~(\ref{5}) in \S3.3, we see that the bracket is zero on $Y_{1,t}$.

\begin{lem}  In the finite case, $\Psi|_{Y_{1,t}}\equiv -\F{(2t-1)}{t_k}$.\end{lem}

{\bf Proof}  The proof that $\Psi_3|_{Y_{1,t}}\equiv 0$ in the affine case applies here without modification.  Lemmas~\ref{psi1lem} and \ref{psi2lem} then apply to give us the result. $\Box$

Comparing to Eq.~(\ref{6}) in \S3.3, we have
$$[v^+\;v^+\;f_{\gamma_1}f_{\gamma_2}.v^+\ldots\;f_{\nu_1}f_{\nu_2}.v^+\;f_{\beta_1}f_{\beta_2}f_{\beta_2}.v^+\;]:=$$$$-\Psi(v^+\otimes [\;v^+\;f_{\gamma_1}f_{\gamma_2}.v^+\ldots\;f_{\nu_1}f_{\nu_2}.v^+\;f_{\beta_1}f_{\beta_2}f_{\beta_2}.v^+\;])-$$$$\F{2t-1}{t_k}[v^+\;v^+\;\;f_{\gamma_1}f_{\gamma_2}.v^+\ldots\;f_{\nu_1}f_{\nu_2}.v^+\;f_{\beta_1}f_{\beta_2}f_{\beta_2}.v^+\;]=0.$$
Note that, by assumption, all the other terms of Eq.~(\ref{6}) are zero: they involve brackets associated to weights that are too high, as in the proof of Lemma~\ref{affpsi}.

\begin{lem}
In the affine case, $\Psi_2|_{U_{1,t}}\equiv (\rho,\alpha)-1$.
\end{lem}

{\bf Proof}  We start by establishing that highest weight vectors in $U_{1,t}$ can take a certain form.  Choose $\alpha\in\Deltak_+$ with minimal height so that $f_\alpha.v^+$ and $f_\alpha.v^+_t$ are both nonzero, $v^+_t=[v^+\;u]$ a highest weight vector of $\g_t$.  Let $w_Y\in Y_{1,t}$ and $w_U\in U_{1,t}$ satisfy
$$f_\alpha.v^+\otimes [v^+\;w]=w_Y+w_U.$$
By weight considerations, $w_U$ is a highest weight vector of $U_{1,t}$, thus generates an irreducible submodule of $U_{1,t}$.  We assume that $U_{1,t}$ is itself irreducible as doing so does not change our argument materially.
We have
$$\Psi_1(f_\alpha.v^+\otimes v^+_t)=\sum_i X_i.f_\alpha.v^+\otimes [Y_i.v^+\;w]+ X_i.f_\alpha.v^+\otimes [v^+\;Y_i.w].$$
Consider terms of the form
\begin{equation}\label{16}
X_i.f_\alpha. v^+\otimes [Y_i.v^+\;w].\end{equation}
Such an expression is nonzero only when $X_i$ is from $\hk$ or a positive root space: when $X_i$ is from a negative root space, $Y_i.v^+=0$.  If $\{H_i\}$ and $\{H'_i\}$ are dual bases of $\hk$ we have
  $$\sum_i H_i.f_\alpha.v^+\otimes [H'_i.v^+\;w]=(\lambda_1,\lambda_1-\alpha)f_\alpha.v^+\otimes [v^+\;w]=f_\alpha.v^+\otimes [v^+\;w].$$ 

Recall next that there is only one simple root $\alpha_i$ (per simple component of $\gk$) that satisfies $f_{\alpha_i}.v^+\neq 0$.  Moreover, the root $\alpha$ appearing in (\ref{16}) has minimal height with  $f_\alpha.v^+$ and $f_\alpha.[v^+\;w]$ both nonzero.  In particular, $\alpha_i$ must be in the support of $\alpha$, but with multiplicity one. 

If $X_i$ is from a positive root space, $X_i.f_\alpha.v^+$ is zero when the positive root associated to $X_i$ has $\alpha_i$ in support, unless $X_i$ is associated to $\alpha$ itself.  We lose no generality in assuming $f_\alpha$ is among the $Y_i$ basis vectors.  The nonzero terms in (\ref{16}) then reduce to 
$$f_\alpha.v^+\otimes [v^+\;w]+v^+\otimes[f_\alpha.v^+\;w].$$
This gives us
$$\Psi_1(f_\alpha.v^+\otimes [v^+\;w])=f_\alpha.v^+\otimes [v^+\;w]+v^+\otimes[f_\alpha.v^+\;w]+$$$$\sum_i X_i.f_\alpha.v^+\otimes [v^+\;Y_i.w]=(1-(\rho,\alpha))w_U+w_Y.$$
Then
$$-(\rho,\alpha)w_U=v^+\otimes[f_\alpha.v^+\;w]+\sum_iX_i.v^+\otimes [f_\alpha.v^+\;Y_i.w].$$

If we take 
$$v^+\otimes[f_\alpha.v^+\;w]=v_U+v_Y,$$
where $v_U\in U_{1,t}$ and $v_Y\in Y_{1,t}$, we can apply $\Psi_1$ to show that 
$$-(\rho,\alpha)v_U=f_\alpha.v^+\otimes [v^+\;w]+\sum_i X_i.v^+\otimes [f_\alpha.v^+\;Y_i.w].$$ 

Let $\kappa$ be the scalar value of $\Psi_2$ on $U_{1,t}$.  We have
$$\Psi_2(w_U+w_Y)=\kappa w_U-w_Y=\sum_i X_i.v^+\otimes [f_\alpha.v^+\;Y_i.w]$$
and
$$\Psi_2(v_U+v_Y)=\kappa v_U-v_Y=\sum_i X_i.f_\alpha.v^+\otimes [v^+\;Y_i.w].$$
Then
$-(\rho,\alpha)v_U=(\kappa+1)w_U$
and
$-(\rho,\alpha)w_U=(\kappa+1)v_U$
imply that $\kappa+1=\pm (\rho,\alpha)$, or, equivalently, that $w_U=\mp v_U$.  Suppose $w_U=v_U$.  We then have 
$$z=w_U+w_Y-(v_U+v_Y)=f_\alpha.v^+\otimes [v^+\;w]-v^+\otimes[f_\alpha.v^+\;w]\in Y_{1,t}.$$
Note that $e_i.[v^+\;w]=[v^+\; e_i.w]=0$ for any positive root vector $e_i\in\gk$, by assumption. Suppose $\alpha$ is simple.  Then $e_i.z=0$ for all $e_i$ so
$z$, a weight vector with weight $\lambda_1+\lambda_t-\alpha$, is a highest weight vector of $Y_{1,t}$.  That is impossible since $Y_{1,t}$ is irreducible with highest weight $\lambda_1+\lambda_t$.  If $\alpha$ is not simple, choose a positive root vector $e_i\in\gk$ so that $e_i.f_\alpha.v^+\neq 0$. Let $[e_i\;f_\alpha]=f_\gamma$ and note that the height of $\gamma$ is less than the height of $\alpha$.  We have 
$$e_i.z=f_\gamma.v^+\otimes [v^+\;w] -v^+\otimes [f_\gamma.v^+\;w]\in Y_{1,t}.$$ 
Repeating the procedure as necessary, we get a highest weight vector in $Y_{1,t}$ associated to a weight strictly less than $\lambda_1+\lambda_t$, a contradiction.  Our conclusion is that $w_U=-v_U$, and $\kappa=(\rho,\alpha)-1$.
$\Box$ 

\begin{lem} In the finite case, $\Psi_2|_{U_{1,t}}\equiv (\rho,\alpha)-\F{t-1}{t_k}-1$.\end{lem}

{\bf Proof} 
As above, we take $f_\alpha.v^+\otimes[v^+\;w]=w_U+w_Y$.  On the one hand, 
we have
$$\Psi_1(f_\alpha.v^+\otimes [v^+\;w])=\sum_i X_i.f_\alpha.v^+\otimes Y_i.[v^+\;w]=$$$$
(\lambda_1-\alpha,\lambda_1)f_\alpha.v^+\otimes [v^+\;w]+v^+\otimes[f_\alpha.v^+\;w]+\sum_i X_i.f_\alpha.v^+\otimes [v^+\;Y_i.w]=$$$$
\left(1-\F{1}{t_k}\right)f_\alpha.v^+\otimes [v^+\;w]+v^+\otimes[f_\alpha.v^+\;w]+\sum_i X_i.f_\alpha.v^+\otimes [v^+\;Y_i.w].$$
On the other hand, we have
$$\Psi_1(f_\alpha.v^+\otimes [v^+\;w])=\Psi_1(w_U)+\Psi_1(w_Y)=\left(1-\F{t}{t_k}-(\rho,\alpha)\right)w_U+\left(1-\F{t}{t_k}\right)w_Y.$$
Subtract 
$$\left(1-\F{t}{t_k}\right)f_\alpha.v^+\otimes[v^+\;w]=\left(1-\F{t}{t_k}\right)(w_U+w_Y)$$
from both expressions to get
\begin{equation}\label{19}-(\rho,\alpha)w_U=\left(\F{t-1}{t_k}\right)f_\alpha.v^+\otimes [v^+\;w]+v^+\otimes [f_\alpha.v^+\;w]+\sum_iX_i.f_\alpha.v^+\otimes [v^+\;Y_i.w].\end{equation}
Applying the same calculations to $w_U+w_Y=v^+\otimes [f_\alpha.v^+\;w]$ we get
\begin{equation}\label{20}-(\rho,\alpha)v_U=\left(\F{t-1}{t_k}\right)v^+\otimes [f_\alpha.v^+\;w]+f_\alpha.v^+\otimes [v^+\;w]+\sum_iX_i.v^+\otimes [f_\alpha.v^+\;Y_i.w].\end{equation}

Let $\tau$ be the scalar value of $\Psi_2$ on $U_{1,t}$, assuming $U_{1,t}$ is irreducible.  We have
$$\Psi_2(w_U+w_Y)=\tau w_U+\left(-1-\F{t-1}{t_k}\right)w_Y=\sum_i X_i.v^+\otimes [f_\alpha.v^+\;Y_i.w],$$
and $$\Psi_2(v_U+v_Y)=\tau v_U+\left(-1-\F{t-1}{t_k}\right)v_Y=\sum_i X_i.f_\alpha.v^+\otimes [v^+\;Y_i.w].$$

Rewrite (\ref{19}) to get
$$-(\rho,\alpha)w_U=\F{t-1}{t_k}(w_U+w_Y)+(v_U+v_Y)+\left(-1-\F{t-1}{t_k}\right) v_Y+\tau v_U.$$
This gives us
$$\left(-(\rho,\alpha)+\F{1-t}{t_k}\right)w_U=(1+\tau)v_U+\F{t-1}{t_k}(w_Y-v_Y).$$
Since $t>1$, $w_Y-v_Y=0$ giving us 
$$\left(-(\rho,\alpha)+\F{1-t}{t_k}\right)w_U=(1+\tau)v_U.$$
The same trick applied to (\ref{20}) gives us
$$\left(-(\rho,\alpha)+\F{1-t}{t_k}\right)v_U=(1+\tau)w_U,$$ which implies
$$\left(-(\rho,\alpha)+\F{1-t}{t_k}\right)^2=(1+\tau)^2.$$
Then $\left(-(\rho,\alpha)+\F{1-t}{t_k}\right)=\pm (1+\tau)$, or equivalently, $w_U=\pm v_U$.  We established that $w_Y=v_Y$, though, so if $w_U=v_U$, we have $f_\alpha.v^+\otimes [v^+\;w]=v^+\otimes [f_\alpha.v^+\;w]$, which is absurd.  We conclude that 
$$\left(-(\rho,\alpha)+\F{1-t}{t_k}\right)=- (1+\tau)$$
giving us $\tau=(\rho,\alpha)+\F{t-1}{t_k}-1$, as desired.
$\Box$

\begin{lem} In both the affine and finite cases, $\Psi_3|_{U_{1,t}}\equiv 0$. \end{lem}

{\bf Proof} 
We do the proof for the affine case; there is no significant difference in the finite case.

We have
$$f_\alpha.v^+\otimes [v^+\;w]=w_Y+w_U=$$\begin{equation}\label{21} w_Y-\F{1}{(\rho,\alpha)}v^+\otimes [f_\alpha.v^+\;w]-\F{1}{(\rho,\alpha)}\sum_i X_i.f_\alpha.v^+\otimes [v^+\;Y_i.w].\end{equation}
This is for a positive root $\alpha$ so that $\lambda_1+\lambda_t-\alpha$ is a highest weight of $U_{1,t}$.

We know that $\Psi_3$ is homothetic on an irreducible component of $U_{1,t}$.  As before, we lose no generality assuming that $U_{1,t}$ is irreducible.  Say that $\Psi_3|_{U_{1,t}}\equiv \tau$.  As per Proposition~\ref{mystery}, take
$$[v^+\;w]=[v^+\;f_{\gamma_1}f_{\gamma_2}.v^+\;\ldots\;f_{\nu_1}f_{\nu_2}.v^+\;f_{\beta_1}f_{\beta_2}f_{\beta_2}.v^+\;].$$
We have
$$\Psi_3(f_\alpha.v^+\otimes [v^+\;w])=$$$$\sum_i X_i.f_{\gamma_1}f_{\gamma_2}.v^+\otimes [f_\alpha.v^+\;v^+\;Y_i.[f_{\xi_1}f_{\xi_2}.v^+\;\ldots\;f_{\nu_1}f_{\nu_2}.v^+\;f_{\beta_1}f_{\beta_2}f_{\beta_2}.v^+\;]\;]+\ldots$$
$$+X_i.f_{\nu_1}f_{\nu_2}.v^+\otimes [f_\alpha.v^+\;v^+\;f_{\gamma_1}f_{\gamma_2}.v^+\;\ldots\;f_{\zeta_1}f_{\zeta_2}.v^+\;Y_i.f_{\beta_1}f_{\beta_2}f_{\beta_2}.v^+\;]$$$$-X_i.f_{\beta_1}f_{\beta_2}f_{\beta_2}.v^+\otimes [f_\alpha.v^+\;v^+\;f_{\gamma_1}f_{\gamma_2}.v^+\;\ldots\;f_{\zeta_1}f_{\zeta_2}.v^+\;Y_i.f_{\nu_1}f_{\nu_2}.v^+\;]$$\begin{equation}\label{22}=-\F{\tau}{(\rho,\alpha)}v^+\otimes [f_\alpha.v^+\;w]-\F{\tau}{(\rho,\alpha)}\sum_i X_i.f_\alpha.v^+\otimes [v^+\;Y_i.w].\end{equation}

If the expression on the right side of the last equal sign in Eq.~(\ref{22}) does not include a nonzero multiple component of $v^+\otimes [f_\alpha.v^+\;w]$, then some $X_i.f_\alpha.v^+\otimes [v^+\;Y_i.w]=-v^+\otimes [f_\alpha.v^+\;w]$.  Without loss of generality, we can assume that $v^+\otimes [v^+\;f_\alpha.w]=-v^+\otimes [f_\alpha.v^+\;w]$.  In this case, $$f_\alpha.(v^+\otimes[v^+\;w])=f_\alpha.v^+\otimes[v^+\;w]\in Y_{1,t}.$$  Comparing to (\ref{21}), we are forced to conclude that
$$v^+\otimes [f_\alpha.v^+\;w]+\sum_i X_i.f_\alpha.v^+\otimes [v^+\;Y_i.w]=0,$$
which is impossible unless $U_{1,t}$ is itself zero.

If $U_{1,t}$ is not zero, the expression on the right side of the equal sign in Eq.~(\ref{22}) includes a nonzero multiple of $v^+\otimes [f_\alpha.v^+\;w]$ but the left side does not, as terms on the left side are all of the form $v_\nu\otimes v_\mu$, $v_\nu \in\g_1$ with $\nu$ no greater than $\lambda_1$ less one positive root of $\gk$.  We conclude that $\tau=0$.
$\Box$ 

\begin{lem}
In the finite case, $\Psi|_{U_{1,t}}\equiv \F{-1}{t_k}$.\end{lem}

{\bf Proof}  
The proof follows the lemmas. $\Box$ 

If we define the bracket in the affine case to be
$\Psi$ and in the finite case to be $\Psi+\F{1}{t_k}$, then it is identically zero on $\cT_t$.  Proposition~\ref{mystery} now follows by a weight argument applied to terms of the form (\ref{mysteq}).  This proves that the algorithms actually produce Lie algebras that enjoy the structure determined by $\Dk$ as extended by $\alpha_k$, which completes the proofs of Theorems~\ref{main1} and \ref{main2}.

\begin{flushleft}
Meighan I. Dillon\\
Mathematics\\
Southern Polytechnic State University\\1100 S. Marietta Pkwy, Marietta, GA 30060\\
{\tt mdillon@spsu.edu}\end{flushleft}

\end{document}